\documentclass[12pt]{amsart}
\title{Rigidity of Schottky sets}
\author{Mario Bonk}
\address{Mario Bonk\\Department of Mathematics\\
University of Michigan\\ 2074 East Hall\\ 530 Church Street\\ Ann
Arbor, MI 48109-1043\\USA} \email{mbonk@umich.edu}
\thanks{M.B.\ was supported by NSF grants DMS-0456940 and
DMS-0244421}

\author{Bruce Kleiner}
\address{Bruce Kleiner\\Mathematics Departement\\
Yale University\\ PO Box 208283\\ New Haven, CT 06520-8283
\\USA} \email{bruce.kleiner@yale.edu}
\thanks{B.K.\ was supported by NSF grants DMS-0505610, and DMS-0701515}

\author{Sergei Merenkov\\ \\ \textit{Dedicated to the memory of Juha Heinonen}}

\address{Sergei Merenkov\\Department of Mathematics\\
University of Illinois\\ 1409 W.~Green Street\\ Urbana, IL 61801
\\USA} \email{merenkov@math.uiuc.edu}
\thanks{S.M.\ was supported by NSF grants DMS-0400636, DMS-0703617,   DMS-0244421,
and DMS-0653439.}

\date{February 17, 2011}

\usepackage{amsmath, amssymb, amsthm,latexsym}
\usepackage{graphicx}
\newcommand\C{{\mathbb C}}
\newcommand\oC{\overline {\mathbb C}}
\newcommand\N{{\mathbb N}}
\newcommand\D{{\mathbb D}}

\newcommand\R{{\mathbb R}}
\renewcommand\H{{\mathbb H}}
\newcommand\Sph{{\mathbb S}}
\newcommand\Ca{{\mathcal C}}
\newcommand\geo{\partial_\infty}
\newcommand\dee{\partial}
\newcommand\Mod{\operatorname{mod}}
\newcommand\dist{\operatorname{dist}}
\newcommand\diam{\operatorname{diam}}
\newcommand\id{\operatorname{id}}
\newcommand\inte{\operatorname{int}}
\newcommand\forr{\quad\text{for} \quad}
\newcommand\co{\colon}
\renewcommand\:{\colon}
\newcommand\sub {\subseteq}
\newcommand\ra {\rightarrow}

\newcommand\ga{\gamma}
\newcommand\la{\lambda}
\newcommand\eps{\epsilon}
\newcommand\OC{{\overline{\mathbb C}}}
\newcommand\no{\noindent} 

\newtheorem{theorem}{Theorem}[section]

\newtheorem{proposition}[theorem]{Proposition}
\newtheorem{corollary}[theorem]{Corollary}

\newtheorem{lemma}[theorem]{Lemma}

\theoremstyle{definition}
\newtheorem{example}[theorem]{Example}

\begin{document}


\abstract{We call a complement of a union 
 of at least three disjoint (round) open balls 
in the unit  sphere $\Sph^n$ a Schottky set. We  prove 
that every quasisymmetric homeomorphism   of a Schottky set of spherical 
 measure zero 
to another Schottky set  
is the restriction of a 
M\"obius transformation on $\Sph^n$.
In the other direction we
show that every  Schottky 
set in $\Sph^2$ of positive measure admits  non-trivial 
quasisymmetric maps to other Schottky sets. 

These results are applied to establish rigidity statements for convex subsets of hyperbolic space that have  totally geodesic boundaries.
}
\endabstract

\maketitle

\section{Introduction}\label{s:Intro}
\no 
Let $\Sph^n$ denote the  $n$-dimensional unit  sphere in $\R^{n+1}$ 
equipped with the restriction of the Euclidean metric. 
A \emph{Schottky set} is a subset of $\Sph^n$  whose complement 
is a union of at least three  disjoint metric (i.e., round) open balls. We impose the requirement that a Schottky set has  at least three distinct open balls   as complementary components to rule out cases that are easy to analyze for the type of problems we consider. 
Each Schottky set is endowed with the induced metric from $\Sph^n$.

Let $f\co X\to
Y$ be a  homeomorphism between two metric spaces $(X, d_X)$ and
$(Y, d_Y)$. The map $f$ is called \emph{$\eta$-quasisymmetric}, where 
 $\eta\co [0,\infty)\to[0,\infty)$ is  a homeomorphism,  
if
$$
\frac{d_Y(f(x),f(y))}{d_Y(f(x),f(z))}\leq
\eta\bigg(\frac{d_X(x,y)}{d_X(x,z)}\bigg)
$$
for every triple of distinct points $x, y, z\in X$. We say that
$f$ is \emph{quasisymmetric} if it is $\eta$-quasisymmetric for
some $\eta$.  For more discussion on  quasisymmetric maps and the related classes of quasi-M\"obius and quasiconformal maps see  Section~\ref{s:qc}.

Every M\"obius transformation on $\Sph^n$
 is a quasisymmetric map and 
sends  Schottky sets to  Schottky sets.  
We say that a Schottky set  $S \sub \Sph^n$ is \emph{rigid} if  this is the 
only way to obtain Schottky sets as  
quasisymmetric images of $S$, i.e., if every quasisymmetric 
map of $S$  onto any other Schottky set $S'\sub \Sph^n$ 
is the  
restriction of a M\"obius transformation.  

In this paper  we consider the problem  of  characterizing rigid Schottky sets.  This is motivated  by some recent investigations on uniformization results  for Sierpi\'nski carpets (see \cite[Ch.~7 and 8] {Bo} for a survey).

The case $n=1$ is trivial. Indeed, according to our definition the Schottky sets in $\Sph^1$  are precisely the closed subsets $S$ of $\Sph^1$ with at least three complementary components. If $S$  contains at least four points, then $S$ is not rigid. To see this note that  every smooth diffeomorphism on $\Sph^1$ that  changes the cross-ratio of four points in $S$ is a quasisymmetric map of $S$ to another Schottky set that does not agree with any M\"obius transformation on $\Sph^1$ restricted to
$S$. Therefore, we can assume $n\ge 2$ in the following.

Our main result is the following sufficient condition for rigidity. 
\begin{theorem}\label{T:Rig1}
Every  Schottky set  in $\Sph^n$, $n\ge 2$, 
 of  spherical  measure zero  is rigid.
\end{theorem}

The proof requires considerable preparation and will be completed in Section~\ref{s:szero}. 

It turns out that in dimension $2$ the condition of vanishing spherical measure is also necessary for the rigidity of a Schottky set.

\begin{theorem}\label{T:Rig2}
A Schottky set  in $\Sph^2$ is rigid 
 if and only if it has spherical measure zero.
\end{theorem}

The proof of the necessity part of this statement uses a rather standard quasiconformal deformation argument (see Section~\ref{s:pos}). It is based on the  measurable Riemann mapping theorem which is  only 
available for $n=2$.

It seems unlikely that a similar simple characterization for the rigidity of a Schottky set can be given in dimensions $n\ge 3$. 
Schottky sets with non-empty interior are always 
non-rigid. It is not hard to construct examples of non-rigid Schottky sets with empty interior in all dimensions (see Example~\ref{ex:nonrig}). By Theorem~\ref{T:Rig1} they necessarily have positive measure.  On the other hand, for $n\ge 3$ there exist rigid Schottky sets $S\sub \Sph^n$ of positive measure. 

\begin{theorem} \label{T:exposrig} 
For each   $n\geq 3$ there exists a  Schottky set in $\Sph^n$ that has positive measure and is rigid.
\end{theorem}

We will construct such sets in Section~\ref{s:posrig}.

Theorem~\ref{T:Rig1}  can be applied to obtain rigidity statements for convex subsets of hyperbolic $n$-space $\H^n$ that have totally geodesic boundary. 

\begin{theorem} \label{T:Rig3}
 Let  $K$ and $K'$ be closed convex sets in $\H^n$, 
$n\ge 3$, with non-empty interior.  Suppose that each set has non-empty 
 boundary consisting of disjoint hyperplanes, and 
that $\partial_\infty K\sub \partial_\infty \H^n\cong \Sph^{n-1} $ has measure zero. 

Then every quasi-isometry between   $K$ and $K'$ has finite distance to the restriction to $K$ of 
 an isometry of $\H^n$ mapping   $K$ to  $K'$.       
\end{theorem} 

In particular, $K$ and $K'$ are isometric. The relation to Theorem~\ref{T:Rig1} is given by the fact that the boundaries at infinity 
 $\geo K$ and $\geo K'$ are Schottky sets and the given quasi-isometry between $K$ and $K'$ induces a quasisymmetric map between $\geo K$ and $\geo K'$ (cf.~Proposition~\ref{prop:qisom}).

The underlying rigidity questions for convex sets in hyperbolic space  will be studied in Section~\ref{s:con}.  We ruled out  $n=2$ in the previous   
theorem, because the statement is not true in this case. For complete results for the low-dimensional cases $n=2$ and $n=3$ 
see  Theorems~\ref{T:Rig4} and~\ref{T:Rig5}. 

Let  $\Gamma$ and $\Gamma'$ be the groups of hyperbolic isometries generated 
by the reflections in the hyperplanes bounding the sets $K$  and $K'$ as in Theorem~\ref{T:Rig3}, respectively.  If we make the additional assumption that the quasi-isometry in the statement is defined on all of $\H^n$ and is  equivariant with respect to 
$\Gamma$ and $\Gamma'$ in a suitable sense, then Theorem~\ref{T:Rig3} can be deduced from 
results by Sullivan (see Theorem~IV and Section~VII in \cite{dS78}).
It is possible to promote every quasi-isometry on the convex set 
 $K$ to a global equivariant map on $\H^n$ by successive reflections
 in an obvious way. It can be shown that the new map is a quasi-isometry on $\H^n$, but there seems to be no simple proof for this fact.

The issue of equivariance turns out to be the main difficulty in the proof of Theorem~\ref{T:Rig1}. In this case one wants to extend
a given quasisymmetric map $f$  between two Schottky sets in $\Sph^n$ to  a quasisymmetric map on $\Sph^n$ that is equivariant with respect to the groups of M\"obius transformations 
generated by the reflections in the ``peripheral spheres" of the Schottky sets (the boundaries of the balls forming the complementary components). 
We study such ``Schottky groups" in Section~\ref{s:sgroups}. The desired equivariant extension  of $f$  is obtained in Proposition~\ref{prop:eqvar2}.  One of the main ingredients in the proof is
the deep  extension theorem for quasiconformal maps due to Tukia and V\"ais\"al\"a \cite{TV82} (cf.~Theorem~\ref{TV-Thm}).

Theorem~\ref{T:Rig3} was already known for hyperbolic convex sets
$K$ and $K'$  with finite inradius,
and a positive lower bound on the separation between boundary components
\cite{kapkleleesch}.
This includes  universal covers  of  
compact hyperbolic $3$-orbifolds with non-empty totally geodesic 
boundaries \cite{KK00}. In this case the Schottky sets $\geo K$
arising as boundaries  are  homeomorphic to a 
Sierpi\'nski
carpet.  The statement and proof in \cite{kapkleleesch}
were inspired by the work of  R.~Schwartz
on nonuniform lattices in the isometry group of $\H^n$ \cite{schwartz}. 
Schwartz' work leads to 
analogous rigidity
statements for subsets $K,K'\sub\H^n$ which are obtained from
$\H^n$ by deleting certain disjoint collections of horoballs.  
His proof involves several steps: showing that boundary components
are preserved by quasi-isometries, that quasi-isometries can be
extended to $\H^n$, and finally, that the boundary homeomorphism
of the extension is conformal almost everywhere.  The proof in \cite{kapkleleesch}
follows the same outline, only each of the steps is simpler
than in the case of horoball complements.
The failure of Theorem~\ref{T:Rig3} when one drops the inradius condition
(which permits the boundary to have positive measure) was also known 
\cite{kapkleleesch}.   
Other results in this direction  
were obtained by Frigerio 
\cite{fri1, fri2} (we thank C.~Leininger for bringing this work to our 
attention).
Also related to this is the  rigidity problem 
 for circle packings \cite{RS} or for conformal maps of circle domains \cite{HS94}.

The outline of the paper is as follows. In Section~\ref{s:sets} we prove some connectivity properties of Schottky sets, and give a topological characterization of the peripheral spheres of a Schottky set. Section~\ref{s:sgroups} discusses properties of the group obtained by successive reflections in the peripheral spheres of a Schottky set. We also recall some facts about Hausdorff convergence of sets.  In Section~\ref{s:qc} we review  quasiconformal and related maps. 
The material in Sections~\ref{s:sets}--\ref{s:qc}
is quite standard. 

We then prove that a quasisymmetric map 
between Schottky sets has an equivariant extension (cf.~Proposition~\ref{prop:eqvar2}). Combined with a differentiation lemma (cf.~Lemma~\ref{lem:conf})
this will give us a proof of Theorem~\ref{T:Rig1} in Section~\ref{s:szero}.  
After some discussion on   Beltrami coefficients, we 
give a proof of  Theorem~\ref{T:Rig2} in 
Section~\ref{s:pos}.  We also discuss an example of a Schottky set with empty interior that is not rigid (cf.~Example~\ref{ex:nonrig}). A rigid Schottky set  
 of positive measure in $\Sph^n$, $n\ge 3$, is constructed in Section~\ref{s:posrig}. The key is a rigidity statement for ``relative" Schottky sets that is of independent interest (Theorem~\ref{T:relSch}). 
 The topic of the final Section~\ref{s:con} is  rigidity statements for convex sets in hyperbolic space with totally geodesic boundaries.  
 
 \medskip\no 
\textbf{Acknowledgement.} The authors  greatly benefitted from  many conversations  with their friend and colleague Juha Heinonen who recently passed away. The authors learned most of the facts about  quasiconformal and related mappings  as presented in Section~\ref{s:qc} from Juha.  The proof of Proposition~\ref{prop:ext} was also inspired by discussions with Juha. We dedicate this  paper to his memory.

\section{Schottky sets} \label{s:sets}

\no 
We first collect some general facts about Schottky sets $S\sub \Sph^n$. We write such a set in the form 
\begin{equation}\label{stdS}
S=\Sph^n\setminus \bigcup_{i\in I}B_i,
\end{equation}
 where the sets   $B_i$, $i\in I$, are pairwise disjoint 
open balls in $\Sph^n$. Here $I=\{1, \dots, l\}$, $l\ge 3$, if $I$ is finite, and $I=\N=\{1,2,3,\dots\}$ if $I$ is infinite. 
The collection of the balls $B_i$, $i\in I$,  is uniquely 
determined by $S$ as it is the set of components of 
$\Sph^n\setminus S$.  We refer to the $(n-1)$-spheres 
$\partial B_i$ as the {\em peripheral spheres} of $S$.  
These sets are topologically distinguished  as 
Proposition~\ref{prop:peri} will  show. 
First  we will discuss some connectedness properties of   Schottky sets. 

\begin{lemma}\label{lem:ballconn}
Let $S\sub\Sph^n$, $n\ge 2$, be a Schottky set, and $B$  an 
open or  closed ball   in $\Sph^n$. Then $S \cap B$ is path-connected. In particular, $S$ is path-connected. 
\end{lemma}

\medskip \no
{\em Proof.} We write $S$ as in \eqref{stdS}. 
If $x,y\in S\cap B$, there exists an arc $\gamma$ of a circle in $\Sph^n$ that connects $x$ and $y$ and is contained in $B$. 
Let $J\sub I$ be the set of indices $i\in I$ for which
$\gamma$ has non-empty intersection with the  ball $B_i$.  For each $i\in J$  there 
exists a maximal subarc $\gamma_i$ of $\gamma$ with $\ga_i\sub \bar B_i$. Since the balls $B_i$ are disjoint, the arcs $\ga_i$, $i\in J$, are pairwise non-overlapping, i.e., no interior point of one arc belongs to any other arc.   Since  the endpoints of $\ga_i$, $i\in J$, are in $\partial B_i$, we can find an 
arc $\tilde\ga_i\sub \partial B_i\cap B$ with the same endpoints as $\ga_i$. We now replace the subarcs 
$\ga_i$, $i\in J$, of $\ga$ by the arcs $\tilde \ga_i$. 
If suitably parametrized, this gives a path $\tilde \ga$
connecting $x$ and $y$ in $S\cap B$. This is clear if $J$ is finite. If   $J$ is infinite, this follows from the fact that  $\diam(\tilde \ga_i)\to 0$ as $i\in J\to \infty$.  The path-connectedness of $S\cap B$ follows. 
\hfill $\Box$ \medskip

Let $(Z,d)$ be a metric space. We denote the open ball of radius $r>0$ centered at $a \in Z$ by $B(a,r)$. The space $Z$  is called $\la$-{\em linearly locally connected}, $\la \ge 1$,  if 
the following two conditions hold:

\smallskip \noindent
($\la$-$LLC_1$):$\quad$
 If $B(a,r)$ is a ball in $Z$ and $x,y\in B(a,r)$, then there
exists a compact connected set $E\sub B(a,\la r)$ containing $x$ and $y$.

\smallskip \noindent
($\la$-$LLC_2$):$\quad$
 If $B(a,r)$ is a ball in $Z$ and $x,y\in Z \setminus B(a,r)$,
then there exists a compact connected set $E\sub Z \setminus B(a,r/\la)$ containing
$x$ and $y$. 
\smallskip

For future reference we record the following immediate consequence of Lemma~\ref{lem:ballconn}. 

\begin{proposition}\label{prop:llc}
Every Schottky set $S\sub \Sph^n$, $n\ge 2$, is 1-linearly locally connected.
\end{proposition}

\no {\em Proof.}  The facts that $S$ is  $1$-$LLC_1$ and $1$-$LLC_2$ follow from Lemma~\ref{lem:ballconn}  applied to the open ball  $B=B(a,r)$ and the closed  ball $B=\Sph^n\setminus B(a,r)$,
respectively, where $B(a,r)$ is as in the $LLC$-conditions. 
  \hfill $\Box$ \medskip

\begin{proposition}\label{prop:peri}
Let $\Sigma$ be a topological $(n-1)$-sphere contained  in a Schottky 
set $S\sub \Sph^n$, $n\ge 2$. Then $S\setminus \Sigma$ is connected 
if and only if $\Sigma$ is a peripheral sphere of $S$. 
\end{proposition}

For a very similar result see \cite[Lem.~2.1]{fri2}.

 \medskip 
\no 
{\em Proof.}  We write $S$ as in \eqref{stdS}. 

 If $\Sigma=\partial B_i$, $i\in I$,  is a peripheral sphere of  $S$, 
then $S\setminus \Sigma$ is connected. 
Indeed, let $B=\Sph^n\setminus \bar B_i$. Then $B$  is an open  ball in $\Sph^n$, and 
Lemma~\ref{lem:ballconn} shows that $S\cap B=S\setminus \Sigma$ is path-connected, and hence connected. 

Conversely, suppose that $\Sigma$ is an embedded $(n-1)$-sphere in $S$ and $S\setminus \Sigma$ is connected. 
By  the Jordan-Brouwer Separation Theorem \cite[Thm.~15, p.~198]{Sp}, 
 the set $\Sph^n \setminus \Sigma$ has 
 two components each of which has $\Sigma$ as its boundary. Since $S\setminus \Sigma$ is connected, it is contained
in one of the components $K$ of $\Sph^n \setminus \Sigma$. Let $K'$ 
be the other  non-empty component of
$\Sph^n\setminus \Sigma$. 
Then $K'\cap S =\emptyset$, and so
$K'$ is covered by the balls $B_i$, $i\in I$. In particular, there exists
one ball $B=B_j$ in this  collection with $K' \cap B\ne \emptyset$. 
Since $B  \cap \Sigma = \emptyset $, it follows that $B\sub K'$. 
Now  $K'$ is connected and $\partial B \cap K'\sub 
S\cap K'=\emptyset$. Hence $B=K'$. This implies that 
$\partial B = \partial K' = \Sigma$, and  so 
$\Sigma $ is a peripheral sphere of $S$. 
\hfill $\Box$ \medskip

\begin{corollary} \label{perpre}
Let $f\co S \ra S'$ be a homeomorphism between  Schottky sets 
$S$ and $S'$ in $\Sph^n$, $n\ge 2$. 
Then $f$ maps  every peripheral sphere of $S$  onto a peripheral sphere of 
$S'$. 
\end{corollary}

\section{Schottky groups} \label{s:sgroups} 

\no
Suppose  $S\sub\Sph^n$ is  a Schottky set in $\Sph^n$,
$n\ge 2$, written as in \eqref{stdS}. 
For each $i\in I$ let $R_i\: \Sph^n\ra \Sph^n$ be the reflection in the 
peripheral sphere $\partial B_i$. 
The subgroup of the group of all M\"obius transformations on $\Sph^n$
generated by the reflections $R_i$, $i\in I$,  is denoted by 
 $\Gamma_S$ and called  the 
 {\em Schottky group associated with} $S$.  
It consists of all M\"obius transformations $U$ of the form 
$U=R_{i_1}\circ \dots \circ R_{i_k}$, where $k\in \N$ and 
$i_1, \dots, i_k\in I$. Since $R_i^2=\id_{\Sph^n}$, where $\id_{\Sph^n}$ is the identity map on $\Sph^n$,  we can assume that in such a representation for $U$ the sequence of indices $i_1, \dots, i_k$  is {\em reduced}, i.e.,  $i_r\ne i_{r+1}$ for $r=1, \dots, k-1$. 

We set 
\begin{equation} \label{Sinfpart}
S_{\infty}=\bigcup_{U\in\Gamma_S}U(S).
\end{equation} 
This set consists of all the copies of the original Schottky set under the transformations in the group $\Gamma_S$. We will later see (cf.\ remark after Lemma~\ref{lem:nested}) that this is a dense subset of $\Sph^n$.  

For $k\in \N$ and a reduced sequence $i_1, \dots, i_k\in I$
 we define open balls 
$$ B_{i_1\dots i_k}: = (R_{i_1}\circ \dots \circ R_{i_{k-1}}) (B_{i_k}). $$ 
Then the following facts are easy to check: 

\begin{itemize}

\smallskip

\smallskip
\item[(i)] $B_{i_1 \dots i_{k}}
\sub B_{i_1 \dots i_{k-1}}$ for all reduced sequences $i_1, \dots, i_k$ and $k>1$,

\smallskip
\item[(ii)]  for fixed $k\in \N$ the balls 
$ B_{i_1\dots i_k}$,  where $i_1, \dots, i_k$   is any reduced 
sequence in $I$,
are pairwise disjoint,

\smallskip
\item[(iii)] $(R_{i_1}\circ \dots \circ R_{i_{k}})(S) 
= {\bar B}_{i_1 \dots i_k}\setminus 
\displaystyle  \bigcup_{i\in I\setminus \{i_{k}\}} B_{i_1 \dots i_{k} i }
$  for all reduced sequences.

\end{itemize}

The last fact shows that $(R_{i_1}\circ \dots \circ R_{i_{k}})(S)$ 
is a Schottky set whose peripheral spheres are $\partial B_{i_1 \dots i_k}$
and $ \partial B_{i_1 \dots i_{k} i }$, $i\in I\setminus \{i_{k}\}$. 
 
 The reflection in $\partial B_{i_1 \dots i_k}$ is given by 
$$ R_{i_1}\circ \dots \circ R_{i_{k-1}}\circ R_{i_{k}}\circ R_{i_{k-1}}
\circ \dots \circ R_{i_1},$$ and hence belongs to $\Gamma_S$. 

The  following result is probably well-known. We give a proof for the sake of 
completeness. 

\begin{proposition}\label{prop:Sgroup} 
The group $\Gamma_S$ is a discrete group of M\"obius transformations
 with a presentation 
given by the generators $R_i$, $i\in I$, and the relations
$R_i^2=\id_{\Sph^n}$, $i\in I$.  
\end{proposition} 

\no 
{\em Proof.} To show that $\Gamma_S$ is discrete
(in the topology of uniform convergence on $\Sph^n$), it is enough to find $\delta>0$
such that 
\begin{equation}\label{discrete}
 \inf_{U\in \Gamma_S\setminus\{\text{id}_{\Sph^n}\}} \biggl(\max_{x\in \Sph^n}
|U(x)-x| \biggr) \ge \delta,
\end{equation} 
i.e., every element in $\Gamma_S$ different from the identity element 
moves a point in $\Sph^n$ by a definite amount.  

To see this, 
consider the  indices  $1,2,3\in I$, and write the corresponding
complementary component of $S$  as $B_l=B(x_l, r_l)$, $l=1,2,3$. 
Then we can take $\delta=\min\{r_1,r_2,r_3\}$ in \eqref{discrete}. 
Indeed, let $U\in \Gamma_S\setminus\{\id_{\Sph^n}\}$ be arbitrary. 
Then there exist $k\in \N$ and a reduced sequence of indices $i_1, \dots, i_k\in I$ 
such that
$$U= R_{i_1}\circ \dots \circ R_{i_k}. $$ 
There is one index $j\in \{1, 2, 3\}$,  
 such that $j\ne i_1$ and $j\ne i_k$. Then 
$$ U(B_j)= B_{i_1\dots i_kj} \sub B_{i_1}.$$ 
Since $B_j\cap B_{i_1}=\emptyset$, this implies that $U(x_j) \not \in B_{j}$
and so $$ |U(x_j)-x_j|\ge r_j \ge \delta$$ as desired. 
Hence $\Gamma_S$ is discrete. 

The same argument also shows that 
$R_{i_1}\circ \dots \circ R_{i_k}\ne \id_{\Sph^n},$
whenever $k\in \N$ and $i_1, \dots, i_k\in I$ is a reduced sequence. 
 Hence $\Gamma_S$ has a presentation
as stated. 
\hfill $\Box$ \medskip

Before we prove the next proposition, we will first review some facts about Hausdorff convergence of sets that will be useful throughout the paper. 
Suppose $X$ is a metric space, and $A,B\sub X$.
 Then their 
{\em Hausdorff distance} $\dist_H(A,B)$ is defined as the infimum of all $\delta
\in (0,\infty]$ such that 
$$ A \sub N_\delta (B) \text{ and } B \sub N_\delta(A). $$ 
Here
 $$N_\delta(M)=\{ x\in X: \dist(x,M)<\delta\}$$
 is the open $\delta$-neigborhood of a set $M\sub X$. Note that the Hausdorff 
 distance between sets is only a ``pseudo"-distance. Namely, it can happen that 
 $\dist_H(A,B)=0$ for sets $A\ne B$. Actually, $\dist_H(A,B)=0$ if and only if 
 $\bar A=\bar B$.

 A  sequence $(A_k)$   of  sets in $X$ is said to {\em (Hausdorff) converge}
  to a set $A\sub X$, written $A_k\to A$,  if 
$$ \dist_H(A_k,A)\to 0 \text{ as } k\to \infty. $$ 
If  $X$ is compact, then every sequence $(A_k)$  of non-empty 
subsets of $X$ {\em subconverges} to a non-empty closed subset $A$ of $X$   
(i.e., the sequence has a convergent subsequence with limit $A$). 

Suppose $A_k\to A$. Then for each $x\in A$ there exists a sequence $(x_k)$ such that $x_k\in A_k$ and $x_k\to x$. 
Conversely, if  for some $x\in X$ there exist a subsequence
$(A_{k_l})$ of $(A_k)$ and corresponding points $x_{k_l}\in A_{k_l}$ with $ x_{k_l} \to x$ as $l\to \infty$, then $x\in \bar A$. 
In particular, this implies  that if $x\in X\setminus \bar A$, then 
$x\in X\setminus A_k$ for large $k$.
   We will use 
these facts repeatedly in the following. 

 The following lemma is straightforward to prove. We leave the details to the reader. 
 
 \begin{lemma}\label{L:ballconv}
Suppose  $(B_k)$ is a sequence of 
closed balls in $\Sph^n$  with  $B_k\to B\sub \Sph^n$, where $B\sub \Sph^n$ is closed. Then $B$ is  a (possibly degenerate) closed ball, and  we have $\partial B_k\to \partial B$.
If $x\in \inte(B)$, then there exists $\delta>0$ such that 
$B(x, \delta) \sub \inte(B_k)$ for large $k$. 
\end{lemma}

Here we denote by $\inte(M)$ the interior of a set $M$. 
We call a closed ball {\em degenerate} if it has radius $0$ and consists of only one point.

The next lemma shows that 
the radii of the balls $B_{i_1\dots i_k}$ as defined above tend to $0$  uniformly as  $k\to
\infty$.

\begin{lemma}\label{lem:ushrink} For every $\delta>0$ 
only finitely many of the balls 
\begin{equation}\label{assctballs}
B_{i_1\dots i_k}, \quad k\in \N \text{ and } i_1, \dots, i_k
\text { a reduced sequence of indices in } I, 
\end{equation}
have diameter $\ge \delta$.\end{lemma} 

\no {\em Proof.} If this is not the case, then there exist
infinitely many of these balls with diameter $\ge \delta$. 
Then we can find a sequence $(D_l)_{l\in \N}$ of distinct 
balls from the collection in \eqref{assctballs}
such that $\bar  D_l$ Hausdorff converges 
to a non-degenerate closed ball $D_\infty$ in $\Sph^n$ 
as $l\to \infty$. Since every ball in \eqref{assctballs}
contains balls of fixed size in its complement (namely one of the balls 
$B_1$ or $B_2$),   we have 
$D_\infty\ne \Sph^n$. Since the boundaries $\partial B_{i_1\dots i_k}$ of the balls in \eqref{assctballs} are distinct sets,  the 
$(n-1)$-spheres  $\Sigma_l=\partial
D_l$ are all distinct. By Lemma~\ref{L:ballconv}  they Hausdorff converge to the $(n-1)$-sphere $\Sigma_\infty:=\partial D_\infty$ as $l\to \infty$. Denote
by $T_l$ for $l\in \N\cup\{\infty\}$  the reflection in the sphere
$\Sigma_l$ on $\Sph^n$. Then $T_l$ converges to $T_\infty$ in the topology
of uniform convergence on $\Sph^n$ as $l\to \infty$. Moreover, 
the reflections  $T_l$, $l\in \N$, are all distinct, 
and  they belong to $\Gamma_S$, because they are 
reflections in spheres bounding balls in 
\eqref{assctballs}. Hence $U_l=T_{l+1}\circ T^{-1}_l
\ne \id_{\Sph^n}$  belongs to $\Gamma_S$ for $l\in \N$, 
and $U_l\to T_\infty\circ T^{-1}_\infty =\id_{\Sph^n} $ as $l\to \infty$.
This contradicts the discreteness of $\Gamma_S$. 
\hfill $\Box$ \medskip

\begin{lemma}\label{lem:nested} 
For each point $x\in \Sph^n\setminus S_\infty$ there exists 
a unique sequence $(i_k)$ in $I$ such that $i_k\ne i_{k+1}$  and 
$x\in B_{i_1\dots i_k}$ for all $k\in \N$. 
\end{lemma} 

 Note that $\diam(B_{i_1\dots i_k})\to 0$ as $k\to 
\infty$ by the previous lemma. Since $\partial B_{i_1\dots i_k}\sub S_\infty$, 
it follows that $S_\infty$ is dense in $\Sph^n$. 

\medskip 
\no {\em Proof.} For existence note that if $x\in \Sph^n\setminus S_\infty$,
then $x\not\in S$. Hence there exists $i_1\in I$ 
such that $x\in B_{i_1}$. Since 
$$x\not \in R_{i_1}(S) = \bar  B_{i_1} \setminus
\bigcup_{i\in I\setminus\{i_1\}} B_{i_1i} \sub S_\infty, $$
there exists $i_2\in I$, $i_2\ne i_1$, such that
$x\in B_{i_1i_2}$. Proceeding in this way, we can inductively define the desired sequence $(i_k)$. 
Uniqueness is clear since for fixed $k\in \N$, the 
balls 
$$ B_{i_1\dots i_k}, \quad i_1, \dots, i_k \text{ is a reduced 
sequence in } I, $$
are pairwise disjoint.
\hfill $\Box$ \medskip

\section{Quasiconformal maps} \label{s:qc}

\no 
We recall
some basic facts about quasiconformal and related mappings (see \cite{Va} for general background).
Let $f\: \Sph^n\ra \Sph^n$ be a homeomorphism, and for $x\in \Sph^n$ 
and   small $r>0$ 
define 
\begin{equation}\label{Lf}
L_f(r,x)=\sup\{|f(y)-f(x)|\co y\in \Sph^n \text{ and } |y-x|=r\}, 
\end{equation}
\begin{equation}\label{lf}
l_f(r,x)=\inf\{|f(y)-f(x)|\co y\in \Sph^n \text{ and } |y-x|=r\}, 
\end{equation}
and the {\em dilatation} of $f$ at $x$ by 
\begin{equation}\label{Hf}
H_f(x)=\limsup_{r\to 0}\frac{L_f(x,r)}{l_f(x,r)}. 
\end{equation}
The map $f$ is called 
{\em  quasiconformal} if
$$
\sup_{x\in\Sph^n}H_f(x)<\infty. 
$$
A quasiconformal map $f$ is called  $H$-{\em quasiconformal}, $H\ge 1$,
if 
$$ H_f(x)\le H \quad \text{
for almost every } x\in \Sph^n. $$ 

Quasiconformality can be defined similarly in other  settings, for example for homeomorphisms between regions in $\Sph^n$ or $\R^n$ or between  Riemannian manifolds.

The composition of an $H$-quasiconformal and an $H'$-quasi\-conformal 
map   is an $(HH')$-quasiconformal map.  
If a homeomorphism  $f$  between regions in $\Sph^n$ is  $1$-quasiconformal, then  $f$ is a conformal (possibly orientation reversing) map. If $n\ge 3$,  then  by Liouville's Theorem $f$ is the restriction of a 
M\"obius transformation (cf.\ \cite[p.~43]{Va}).

If $x_1,x_2,x_3,x_4$ are four  distinct points  
in a metric space $(X,d)$, then   their  {\em cross-ratio} 
is the quantity 
$$[x_1,x_2,x_3,x_4]=\frac{d(x_1,x_3)d(x_2,x_4)}{d(x_1,x_4)d(x_2,x_3)}. $$

Let $\eta\:[0,\infty)\ra[0,\infty)$
be a homeomorphism, 
 and  $f\:X\ra Y$ a homeomorphism between metric spaces $(X,d_X)$ and $(Y,d_Y)$.
The map $f$ is an {\em $\eta$-quasi-M\"obius  map} if
$$[f(x_1),f(x_2),f(x_3),f(x_4)]\leq \eta([x_1,x_2,x_3,x_4]).$$
for every $4$-tuple $(x_1,x_2,x_3,x_4)$ of distinct
points in $X$. 

Note that a M\"obius transformation on $\Sph^n$ preserves cross-ratios of points. As a consequence every pre- or post-composition of an $\eta$-quasi-M\"obius map $f\:\Sph^n\ra \Sph^n$ by a M\"obius transformation is $\eta$-quasi-M\"obius.

We list some interrelations between  the classes of maps we discussed
\cite{Va2}:

\begin{itemize}

\item[(i)] Let $n\ge 2$. 
Then every $H$-quasiconformal map $f\:\Sph^n\ra \Sph^n$
is $\eta$-quasi-M\"obius
with $\eta$ depending only on $n$ and $H$. Conversely, 
every $\eta$-quasi-M\"obius map $f\:\Sph^n\to\Sph^n$ is $H$-quasiconformal with 
$H$ depending only on $\eta$.

\item[(ii)] An $\eta$-quasisymmetric map between metric spaces is $\tilde \eta$-quasi-M\"obius
 with $\tilde \eta$ depending only on $\eta$.
\end{itemize}

Conversely, every quasi-M\"obius map between bounded spaces is quasisymmetric.
This statement cannot be  quantitative in general, as follows from the fact that M\"obius transformations act triply transitive on the Riemann sphere.  If  one imposes a three-point normalization on the   
quasi-M\"obius map, then  a quantitative converse is true; indeed, we have:

\begin{itemize}
\item[(iii)] Let $(X,d_X)$ and $(Y,d_Y)$ be bounded metric spaces,
     $f\: X\ra Y$ an  $\eta$-quasi-M\"obius map, $\la \ge 1$, 
$x_1,x_2,x_3 \in X$. 
Set 
$y_i=f(x_i)$, and suppose that 
 $d_X(x_i, x_j)\ge \diam(X)/\la$ and  
$d_Y(y_i, y_j)\ge \diam(Y)/\la$ for  $i\ne j$.  
Then $f$ is $\tilde  \eta$-quasisymmetric with $\tilde  \eta$
depending only on $\eta$ and $\la$.
\end{itemize}

We consider $\R^n$ as a subspace of $\R^{n+1}$ as usual
by identifying a point $(x_1, \dots, x_n)\in \R^n$ with
$(x_1, \dots, x_n,0)\in \R^{n+1}$. In this way, we can also consider 
$\Sph^{n-1}=\Sph^n\cap\R^n$ as a subspace of $\Sph^n$. 

We need the following  deep result 
due to  Tukia and V\"ais\"al\"a~\cite{TV82}.  

\begin{theorem}\label{TV-Thm} Let $n \ge 3$. 
Every $H$-quasiconformal map  $f\: \R^{n-1} \ra \R^{n-1}$   
 has an $H'$-quasiconformal  extension $F\co \R^n \ra \R^n $,
where $H'$ only depends on $n$ and  $H$.
\end{theorem}

For $n=2$ we have  the classical Ahlfors-Beurling extension theorem
that can be formulated as follows. 

\begin{theorem}\label{TV-Thm2}
Every $\eta$-quasisymmetric   map 
$f\co \R\ra \R$ 
 has an $H$-quasiconformal  extension $F\co \R^2 \ra \R^2$,  
where $H$ only depends on $\eta$.    
\end{theorem}

We need the following consequence of these results. 

\begin{proposition} \label{prop:ballext}
 Let $D\ne \Sph^n$ and $D'\ne \Sph^n$ be closed non-degenerate balls in $\Sph^n$, $n\ge 2$,
and $f\:\partial D\ra \partial D'$  a homeomorphism.

\begin{itemize} 

\item[\textrm(i)]
If $f$ is $\eta$-quasi-M\"obius, then it 
can be extended to 
an $\eta'$-quasi-M\"obius  map $F\: D \ra D'$,
where $\eta'$ only depends on $n$ and $\eta$.

\smallskip 
\item[\textrm(ii)]
If each of the balls $D$ and $D'$ 
 is contained in a hemisphere, and $f$ is  $\eta$-quasisymmetric, 
then $f$ can be extended to an   $\eta'$-quasisymmetric map 
 $F\: D \ra D'$,
where $\eta'$ only depends on $n$ and $\eta$. 

\end{itemize} 
\end{proposition} 

\no
{\em Proof.} To prove (i), we  
 map 
$D$ and $D'$ to closed hemispheres by auxiliary M\"obius transformations. 
We may assume that these hemispheres are bounded by 
$\Sph^{n-1}=\Sph^n \cap \R^n \sub \Sph^n$. 
So after suitable composition  of $f$  by M\"obius transformations, we obtain an
$\eta$-quasi-M\"obius map $\tilde f\: \Sph^{n-1}\ra \Sph^{n-1}$. 

If we distinguish suitable points as points at infinity in the two copies of 
 $\Sph^{n-1}$ and make the identification
$\Sph^{n-1}=\R^{n-1}\cup \{\infty\}$, then $\tilde f(\infty)=\infty$, 
and $\tilde f$ restricts to an $\eta$-quasi-M\"obius map 
$ \tilde f\: \R^{n-1} \ra  \R^{n-1}$.   Here $\R^{n-1}$ 
has to be considered as equipped with the chordal metric coming 
from the identification of $\Sph^{n-1}$ with $\R^{n-1}\cup \{\infty\}$
by stereographic projection. 
Cross-ratios for points in $\R^{n-1}$ are  the same 
if we take the chordal metric or the Euclidean metric.  
It follows that $\tilde f : \R^{n-1} \ra \R^{n-1}$ 
is $\eta$-quasi-M\"obius if $\R^{n-1}$ is equipped with the 
Euclidean metric.  Since $\tilde f(\infty)=\infty$, we conclude by
a limiting argument that $\tilde f\: \R^{n-1}\ra \R^{n-1}$ 
is also  $\eta$-quasisymmetric when $\R^{n-1}$ carries this metric.

If $n\ge 3$, this implies that $\tilde f\:\R^{n-1}\ra \R^{n-1}$
is $H$-quasiconformal 
with $H$ only depending on $\eta$. Hence by Theorem~\ref{TV-Thm},  
$\tilde f$ has an $H'$-quasiconformal extension $\tilde F\:\R^n \ra \R^n$ 
with $H'$ depending only on $n$ and $H$, and hence only on $n$ and 
$\eta$. 
If $n=2$, then we get such an  $H'$-quasiconformal extension $\tilde F$ 
from   the Ahlfors-Beurling  Theorem~\ref{TV-Thm2}.  

Letting $\tilde F(\infty)=\infty$ and making the identification 
$\Sph^n=\R^n \cup\{\infty\}$, we get an $H'$-quasiconformal 
mapping $\tilde F\: \Sph^n\ra \Sph^n$ that extends $\tilde f\: \Sph^{n-1}
\ra \Sph^{n-1}$. Note that points are ``removable singularities" for quasiconformal maps \cite[Thm.~17.3]{Va}. Moreover, the dilatation of $\tilde F$ 
 does not change by the passage 
from the Euclidean metric on $\R^n$ to the chordal metric on $\R^n  \sub
\Sph^n=\R^n\cup\{\infty\}$, because these metrics are ``asymptotically" conformal.    
Then $\tilde F$ will be $\eta'$-quasi-M\"obius with $\eta'$ only depending 
on $n$ and $H'$ and hence only on $n$ and $\eta$. Conjugating this extension 
back by the auxiliary M\"obius transformations used above, and 
restricting to a map on $D$, we get an extension of $f$ with the desired
properties.  

To prove  part (ii) suppose that $f$ is $\eta$-quasisymmetric.
Since quasisymmetric maps 
are quasi-M\"obius maps quantitatively, it follows from the first part 
of the proof that there exists
an $\tilde \eta$-quasi-M\"obius extension $F\: D \ra D'$, 
where  $\tilde \eta$ only depends on $n$ and $\eta$. 
If  $D$ and $D'$ are each  contained in a hemisphere, then 
$\diam(D)=\diam(\partial D)$ and $ \diam(D')=\diam(\partial D')$. 
Pick points $x_1,x_2,x_3 \in \partial D$ such that 
$$ |x_i-x_j|\ge \diam(\partial D)/2 =\diam(D)/2  \forr i\ne j,$$ 
and define $y_i=F(x_i)=f(x_i)\in \partial D'$.  
Now by the quasisymmetry of $f$, 
$$ |f(z)-f(x_i)|\le \eta(2) |f(x_j)-f(x_i)|$$ 
for arbitrary $i\ne j$ and $z\in \partial D$.  
It follows that
$$  \diam(D')/\lambda=\diam(\partial D')/\lambda \le |y_i-y_j| \forr i\ne j,$$
where $\la = 2\eta(2)$. Since $\lambda$ only depends on $\eta$, 
it follows from fact (iii) above that  $F$ is $\eta'$-quasisymmetric with $\eta'$ only depending 
on $n$ and $\eta$. 
\hfill $\Box$ \medskip

\section{Extension of quasisymmetric maps between Schottky sets} 
\label{s:ext}

\no 
Throughout this section $S$ and $S'$ will be Schottky  sets in $\Sph^n$,
$n\ge 2$,  
such that there exists a quasisymmetric map
$f\co S\to S'$. By Corollary~\ref{perpre} 
 we can write
\begin{equation} \label{normform} 
S=\Sph^n\setminus \bigcup_{i\in I}B_i \quad\text{and}\quad
S'=\Sph^n\setminus \bigcup_{i\in I}B'_i,  
\end{equation} 
 where both  collections $\{B_i: i\in I\}$ and $\{B'_i: i\in I\}$
consist of pairwise disjoint
open balls in $\Sph^n$, and $f(\partial B_i)=\partial B'_i$ for $i\in I$.
For $i\in I$ let $R_i$ be the reflection in $\partial B_i$, and $R_i'$ 
be the reflection in $\partial B'_i$. 
If $U$ is an element in the Schottky group $\Gamma_S$, then it can be uniquely 
written as
$$ U= R_{i_1}\circ \dots \circ R_{i_k}, $$ 
where $k\in \N$ and $i_1, \dots, i_k$ is a reduced sequence in $I$.

By Proposition~\ref{prop:Sgroup}  the map $\Phi\: \Gamma_S \ra \Gamma_{S'}$ given by 
$$\Phi(U)=U':= R'_{i_1}\circ \dots \circ R'_{i_k}$$ 
is well-defined and defines a group isomorphism between $\Gamma_S$ and $\Gamma_{S'}$. 

Let 
$$ 
S_\infty= \bigcup_{U\in \Gamma_S} U(S) \quad \text{and}  \quad
S'_\infty= \bigcup_{V\in \Gamma_{S'}} V(S'). $$ 
Then $S_\infty$ and $S'_\infty$ are dense sets in $\Sph^n$, 
  the set $S_\infty$ is invariant  under the group $\Gamma_S$, 
 and $S'_\infty$  under $\Gamma_{S'}$. 

\begin{lemma} \label{lem:eqvar1} 
There exists a unique  bijection $f_\infty \:S_\infty \ra  S'_\infty$ 
that extends $f$ equivariantly, that is, $f_\infty|S =f$ and 
$f_\infty\circ U = U'\circ f_\infty$ for all $U\in \Gamma_S$. 
\end{lemma} 

\no {\em Proof.} 
 Let $z\in S_\infty$ be arbitrary. Then there exist $x\in S$ and $U\in \Gamma_S$ such that $z=U(x)$. We define 
  $f_\infty\: S_\infty \ra S_\infty'$ by setting 
  $f_\infty(z)=U'(f(x))$. 

To show that $f_\infty$ is well-defined 
assume that
$U(x)=V(y)$, where $x,y\in S$ and $U,V\in \Gamma_S$, $U\ne V$.  
Then $U^{-1}\circ V=R_{i_1}\circ \dots\circ  R_{i_k}$ and hence 
  $x=R_{i_1}\circ \dots \circ R_{i_k}(y)$,  where $k\in \N$ and $i_1,\dots, i_k$ is a reduced sequence in $I$. It follows that $x\in \bar B_{i_1}\cap S=
  \partial  B_{i_1}$. Thus $x=R_{i_1}(x),$ and so 
  $x=R_{i_2}\circ \dots \circ R_{i_k}(y)$. Repeating this argument, we 
  deduce that $x$ lies on all the spheres $\partial B_{i_1}, \dots ,\partial 
  B_{i_k}$, 
  and is fixed by each of the reflections $R_{i_1},  \dots ,R_{i_k}$. 
  This  shows that $x=y$. Therefore, $f(x)=f(y)$ is fixed by each 
  of the reflections  $R'_{i_1}, \dots, R'_{i_k}$, and so 
  $f(x)=R'_{i_1}\circ \dots \circ R'_{i_k} (f(y))$. 
  Since $R'_{i_1}\circ \dots \circ R'_{i_k}=U'\circ V'^{-1}$, we conclude 
  $U'(f(x))=V'(f(y))$. This implies that $f_\infty$ is well-defined. 
  
It is clear that $f_\infty$ is the unique equivariant extension of $f$ to $S_\infty$. 
An inverse map for $f_\infty$ can be defined similarly.    So 
$f_\infty$ is indeed  a bijection.  
\hfill $\Box$ 
\medskip

The argument in the previous proof also shows that if two copies $U(S)$ and $V(S)$, 
$U,V\in \Gamma_S$, $U\ne V$, of the Schottky set $S$  have a common point $z$, then $z$ lies on peripheral spheres of $U(S)$ and  $V(S)$. 
Note that in general these peripheral spheres need not be identical, but they can be distinct spheres that touch 
at $z$. 
In any case, $U(S)$ and $V(S)$ intersect in a set of measure zero. 
Therefore,  the representation of $S_\infty$ as in \eqref{Sinfpart} 
gives a measurable partition of this set. This will be important 
in the proof of Theorem~\ref{T:Rig2}.  

We would like to prove that $f$ actually has an equivariant quasiconformal
extension to $\Sph^n$. This would easily follow from the   previous lemma
if we could show that $f_\infty$  is a quasi-M\"obius map. 
Though this is true, there seems to be  no straightforward 
proof of this fact. 

We will address this issue by first extending $f$ in a non-equivariant 
way to a quasiconformal map on $\Sph^n$, 
and then correcting this map successively to make it equivariant while 
keeping a uniform bound on the dilatation of the intermediate 
 quasiconformal maps.  The equivariant extension is then obtained as a 
 sublimit of these maps.

The first step is provided by the following 
extension result.
\begin{proposition}\label{prop:ext}
Every quasisymmetric map between  Schottky sets  in $\Sph^n$, $n\ge 2$,
 extends to a quasiconformal homeomorphism  of  $\Sph^n$.
\end{proposition}

\no \emph{Proof.} Suppose  
$f\co S\to S'$ is  an $\eta$-quasisymmetric map between two 
Schottky sets $S$ and $S'$ in $\Sph^n$, $n\ge 2$. We can write 
$S$ and $S'$ as in \eqref{normform}. 
Moreover, by applying suitable M\"obius transformations to $S$ and $S'$ 
if necessary, we may assume that each of the balls $B_i$ and 
$B'_i$, $i\in I$, is contained in a hemisphere. 

By Proposition~\ref{prop:ballext} we can  
extend each map $f|\partial B_i
\co \partial B_i \ra \partial B'_i$, $i\in I$, to 
an  $\eta'$-quasisymmetric map  of  $\bar  B_i$ onto   
$\bar  B'_i$, where $\eta'$ is independent of $i$.
These maps paste together to  a bijection $ F\co \Sph^n \ra \Sph^n$ 
whose restriction to $S$ agrees with $f$ and whose restriction 
to each ball $\bar B_i$ is an $\eta'$-quasisymmetric  map onto $
\bar B_i'$.

We claim that this  global map $F$  
is a quasiconformal homeomorphism.   
We need to show $F$ is continuous and that there exists a constant $H\ge 1$ such that for every
$x\in\Sph^n$, 
\begin{equation}\label{E:Qc}
\limsup_{r\to 0}\frac{L_F(x,r)}{l_F(x,r)}\leq H,  
\end{equation}
where $L_F$ and $l_F$ are defined as in \eqref{Lf} and \eqref{lf}. 
Below  we will write $a\lesssim b$ for two quantities $a$ and $b$ 
if there exists
a constant $C$ that depends only on
the functions $\eta$ and $\eta'$, such that $a\leq Cb$. We will write $a\simeq b$
if both $a\lesssim b$ and $b\lesssim a$ hold.

If $x$ is in a complementary component of $S$, 
then $F$ is continuous at $x$ and (\ref{E:Qc})
follows   from the definition of $F$ with $H=\eta'(1)$. Thus it is enough 
to consider only the case $x \in S$.

Since $S$ is 
connected, there exists small $r_0>0$ such that the spheres 
$$\Sigma(x,r):= \{y\in \Sph^n\co |y-x|=r\}$$
intersect $S$ for $0<r\leq r_0$. Suppose that $r$ is 
in this range and $y\in \Sigma(x,r)$ is  arbitrary. 
Since $F|S=f$ is 
$\eta$-quasisymmetric, it suffices  to show that there exist
points $v', v'' \in S \cap \Sigma(x,r)$   
such that 
\begin{equation}\label{E:Ineq}
|F(v'')-F(x)|\lesssim|F(y)-F(x)|\lesssim|F(v')-F(x)|.
\end{equation}
For then the continuity of $F$ at $x$ will follow from the continuity of $f=F|S$ at $x$, and 
$L_F(x,r)/l_F(x,r)$ will be bounded by a quantity 
 comparable to $\eta(1)$. 

This is trivial if $y$ itself is in $S$. Thus we assume that $y$ is not 
in $S$, i.e., it lies in one of the complementary components of $S$, which
we denote by $B$. 
Let $v'$ denotes an arbitrary point which is in the intersection
of the sphere $\Sigma (x,r)$ and $\partial B$, and 
let $u'$ be the point in the intersection of $\dee B$ and
the geodesic  segment $[x,y]$ (with respect to the spherical metric). Since
 $|y-u'|\leq|v'-u'|,\
|u'-x|\leq|v'-x|$ and $|v'-u'| \le 2r=2|v'-x|$, the triple
$\{x, v', u'\}$ is in $S$, and the triple $\{y, v', u'\}$ is in 
$\bar B$, we have
\begin{align}
|F(y)-F(x)|&\leq|F(y)-F(u')|+|F(u')-F(x)| \notag\\ 
&\lesssim|F(v')-F(u')|+|F(v')-F(x)|\lesssim|F(v')-F(x)|.\notag
\end{align}
This shows the right-hand side of~(\ref{E:Ineq}). To prove the left-hand
side of the inequality, we choose $v''$ in the same way as $v'$, namely to be an 
arbitrary point in the intersection of the sphere $\Sigma(x,r)$ and $\dee B$.
As for $u''$, we choose it to be the preimage under $F$ of  
the intersection point of the geodesic segment $[F(x),F(y)]$ with  $F(\dee B)$.
Again, the triple
$\{x, v'', u''\}$ is in $S$, and the triple $\{y, v'', u''\}$ 
is in  $\bar B$.
We need to consider two cases:
\newline
\emph{Case 1.} $|u''-x|\geq\frac12 r$. In this case we have 
$|v''-x|=r\le 2|u''-x|$, and therefore
$$
|F(v'')-F(x)|\lesssim|F(u'')-F(x)|\leq|F(y)-F(x)|.
$$
\newline
\emph{Case 2.} $|u''-x|\leq\frac12 r$. Then we have 
$|y-u''|\ge \frac12 r$. So 
$$|v''-u''|\le |v''-x|+|u''-x|\le \frac 32 r \le 3 |y-u''|,$$ and thus
\begin{align}
|F(v'')-F(x)|&\leq|F(v'')-F(u'')|+|F(u'')-F(x)|\notag\\
&\lesssim|F(y)-F(u'')|
+|F(u'')-F(x)|=|F(y)-F(x)|.\notag
\end{align}
This completes the proof of~(\ref{E:Ineq}), and thus of~(\ref{E:Qc})
and the proposition. $\Box$
\medskip

Suppose $T \sub \Sph^n$, $n\ge 2$, is a Schottky set, $\Sigma$ one of the 
peripheral spheres of $T$, and $R$ the reflection in $\Sigma$. 
Then $\tilde T= T \cup R(T)$ is also a Schottky set, called the {\em double}
of $T$ along $\Sigma$. 
Let $T'$ be another Schottky set in $\Sph^n$, and $F\: \Sph^n \ra \Sph^n$
be an $H$-quasiconformal map with $F(T)=T'$. Then $\Sigma'=F(\Sigma)$ 
is a peripheral sphere of $T'$. Let $R'$ be the reflection 
in $\Sigma'$, and  $\tilde T'$ be the double 
of $T'$ along $\Sigma'$. Denote by $B$ the open 
ball with $\Sigma=\partial B$ and $B\cap T=\emptyset$. We 
define a map $\tilde F\: \Sph^n \ra \Sph^n$  by 
$$ \tilde F(x)=
\left\{
\begin{array} {ll}
F(x),                  &  \mbox {$x\in \Sph^n\setminus B$}, \\
R'\circ F \circ R(x),  &  \mbox {$x\in \bar B$. } 
\end{array}
\right. 
$$ 
Note that this definition is consistent on $\partial B=\Sigma$ and hence
defines a homeomorphism  from $\Sph^n$ onto itself. 

\begin{lemma} \label{lem:doub} 
The map $\tilde F$ is  an $H$-quasiconformal map 
with $\tilde F|T = F$,  $\tilde F(\tilde T)= \tilde T'$, and $\tilde F \circ 
R= R'\circ \tilde F$. 
\end{lemma}

The main point here is that we get the same dilatation bound for $\tilde F$ as
for $F$. 
In other words, if there exists an $H$-quasiconformal homeomorphism 
of $\Sph^n$ mapping a Schottky set $T$ to a Schottky set $T'$, 
then there also exists a natural 
$H$-quasiconformal homeomorphism that takes a double of $T$ to the corresponding
double of $T'$ and agrees with the original map on $T$.

\medskip 
\no 
{\em Proof.} 
Since M\"obius transformations  are $1$-quasiconformal, the map  $R'\circ F
\circ R$ is $H$-quasiconformal. Hence $\tilde F | B$ and
$\tilde F | (\Sph^n\setminus \bar B)$ are $H$-quasiconformal. 
This implies that $\tilde F$ is $H$-quasiconformal, because sets of 
$\sigma$-finite Hausdorff $(n-1)$-measure (such as $\Sigma=\partial B$) form 
``removable singularities" for quasiconformal maps on $\Sph^n$ (see \cite[Sect.~35]{Va}). 
The other statements are obvious.  
\hfill $\Box$ \medskip

 With the setup as in the beginning of the section we can now prove the following lemma. 

\begin{lemma}  \label{lem:succdoub} 
There exist  $H\ge 1$, Schottky sets $S_k$ and $S'_k$ in $\Sph^n$, and
$H$-quasiconformal  maps $F_k \: \Sph^n \ra \Sph^n$ for $k\in \N_0$ with the following properties: 

\begin{itemize}

\smallskip 
\item[\textrm{(i)}]  $F_0=F$, $S_0=S$, $S_0'=S'$,   

\smallskip 
\item[\textrm{(ii)}] $F_k(S_k)=S'_k$ 
for $k\in \N_0$,

\smallskip 
\item[\textrm{(iii)}] 
$S_{k+1} \supseteq S_k$ 
is a double of $S_k$, and $S'_{k+1}\supseteq S'_k$
is the corresponding double 
of $S'_{k}$ for $k\in \N_0$, 

\smallskip 
\item[\textrm{(iv)}] $F_k|S_k=f_\infty | S_k$ for $k\in \N_0$, 

\smallskip
\item[\textrm{(v)}]  $\displaystyle \bigcup_{k\in \N_0}  S_k= S_\infty$.  

\end{itemize}

\end{lemma}

\no {\em Proof.} 
Define $S_0:=S$, $S'_0:=S'$, and let $F_0\: \Sph^n \ra 
\Sph^n$ be a quasiconformal extension of $f$  as provided by 
Proposition~\ref{prop:ext}. The map $F_0$ will be $H$-quasiconformal for some $H\ge 1$.  

Now if Schottky sets $S_k$ and $S'_k$ and an $H$-quasiconformal map $F_k$ on $\Sph^n$ with $F_k(S_k)=S'_k$ have been defined for some 
$k\in \N_0$, 
we let $S_{k+1} $ be the double of $S_k$ along a peripheral sphere
$\Sigma$
 of $S_k$  with the largest radius (which exists, because there are  only finitely many peripheral spheres whose radii exceed a given positive number). Then $S'_{k+1}$ is defined as the double of $S'_k$  along the peripheral sphere that corresponds to $\Sigma$ 
 under $F_k$, and $F_{k+1}$ is the $H$-quasiconformal 
map obtained from $F_k$ and these doubles as in Lemma ~\ref{lem:doub}.

With these definitions the asserted properties (i)--(iii) are clear.
Since $S_k$ is obtained by successive doubles of $S$, every peripheral sphere of $S_k$ is  an image of a peripheral
sphere of $S$ under a M\"obius transformation in $\Gamma_S$.
In particular, any reflection in a peripheral sphere of $S_k$ belongs
to $\Gamma_S$. Using this,  
property (iv)  follows from the definition of $F_{k}$ and the equivariance 
of $f_\infty$ (cf.~Lemma~\ref{lem:eqvar1})  by  induction on $k$.  

Note that if $r_k$ 
is the maximal radius of a peripheral sphere of $S_k$, then 
\begin{equation}\label{maxsph}
r_k \to 0 \text{ as } k\to \infty.
\end{equation}
Indeed, the peripheral spheres of $S_k$ are among the spheres
\begin{equation}\label{perisph}
\partial B_{i_1\dots i_l}=R_{i_1}\circ \dots \circ R_{i_{l-1}}(\partial B_{i_l}),
\end{equation}
where $l\in \N$ 
and $i_1,\dots, i_l$ is a reduced sequence in $I$.  By 
Lemma~\ref{lem:ushrink} there are only finitely many among the spheres
in \eqref{perisph} 
whose radii exceed any given positive constant $\epsilon>0$. 
Since in the construction of $S_{k+1}$ from $S_k$
we double $S_k$ along a peripheral sphere of maximal radius
and this sphere will not be a peripheral sphere of any of the Schottky sets
$S_{k+1}, S_{k+2}, \dots$,  all spheres of radius $\ge\epsilon$ in 
\eqref{perisph} are eventually
eliminated as possible peripheral spheres of the sets $S_k$ in the doubling process. Therefore, \eqref{maxsph} follows.

Now we can show that (v) holds. 
It is clear that $S_\infty$ contains each $S_k$. 
Suppose $\tilde S:=\bigcup_{k\in \N_0}S_k$ is a proper subset of $S_\infty$. Then $\tilde S$ does not contain all the copies of $S$ under the transformations in $\Gamma_S$, and so there  exist $U\in \Gamma_S$ such that $U(S)$ is not contained in $\tilde S$.
Each such map $U$ has a unique representation 
in the form $U=R_{i_1}\circ \dots\circ R_{i_l}$, where $l\in \N$,
and $i_1, \dots, i_l$ is a reduced finite sequence in $I$.
We fix  $U$ so that it  has a representation of this form with minimal $l$ among all such group elements. Then 
$T=R_{i_1}\circ \dots \circ R_{i_{l-1}}(S)$ is a subset of $S_k$ for sufficiently
large $k$, but $T'=U(S)$ is not. Since the Schottky sets $T$ and $T'$ have the common peripheral sphere $\Sigma= 
R_{i_1}\circ \dots \circ R_{i_{l-1}}(\partial B_{i_l})$, 
this is a peripheral sphere of $S_k$ for all sufficiently large $k$. 
This is impossible, because the radius  of the largest peripheral sphere of $S_k$ tends to $0$ as $k\to \infty$.  
\hfill $\Box$ \medskip

\begin{proposition} \label{prop:eqvar2} 
The quasisymmetric map $f\: S\ra S'$ has 
an equivariant quasiconformal extension 
$F\: \Sph^n \ra \Sph^n$, that is,  $F|S=f$ and $F\circ  U = U' \circ 
F$ for all $U\in \Gamma_S$. 
\end{proposition} 

\no {\em Proof.} 
Consider the Schottky sets $S_k$ and $S'_k$ 
and the   $H$-quasiconformal maps $F_k$ obtained in Lemma~\ref{lem:succdoub}.
Since these maps are uniformly quasiconformal,  
there exists a distortion function $\eta$ such that $F_k$ is an
$\eta$-quasi-M\"obius  map for all $k\in \N_0$. 
Any four points  in $S_\infty$ are contained 
in one of the Schottky sets 
$S_k$, $k\in \N_0$. Since $F_k|S_k=f_\infty|S_k$, it follows that
$f_\infty$ is an $\eta$-quasi-M\"obius map from $S_\infty$ to $S'_\infty$.
Since $S_\infty$ and $S'_\infty$ are dense in $\Sph^n$, the 
map  $f_\infty$ has a unique quasi-M\"obius extension 
$F\:\Sph^n \ra \Sph^n$. Then $F$ is a quasiconformal extension of $f$. 
The  map $F$  has the desired equivariance property as follows from the 
corresponding property of $f_\infty$. 
\hfill $\Box$ \medskip

\section{Schottky sets of measure zero} 
\label{s:szero} 

\no Now we are almost ready to prove Theorem~\ref{T:Rig1}. We need one final ingredient. 

\begin{lemma}\label{lem:conf}  Let $g\: \R^n \ra  \R^n$, $n\in \N$,  
be a map that is differentiable  at $0$. 
Suppose there exists a  sequence $(D_k)$ of 
  non-degenerate  closed balls in $\R^n$ with
$\diam (D_k)\to 0$ such that 
$0\in D_k$ and $D'_k=g(D_k)$ is a ball  for all $k\in \N$.  

Then the derivative $Dg(0)$ of $g$ at $0$
is a (possibly degenerate or orientation reversing) conformal linear 
map, i.e., $Dg(0)=\la T$, where $\la \ge 0$ and $T\: \R^n \ra \R^n$ is a linear 
isometry. 
\end{lemma} 

\no {\em Proof.} We may assume that $g(0)=0$. 
Let $r_k>0$ be the radius of $D_k$, and define 
$\tilde D_k= \frac1{r_k} D_k$
for $k\in \N$.  Then $\tilde D_k$
is a  closed ball of radius $1$ containing $0$. 
By passing to a subsequence if necessary, we may 
assume that the balls $\tilde D_k$ 
Hausdorff converge to a closed ball $D\sub \R^n$ 
of radius $1$.  

Since $r_k\to 0$,  the maps 
$g_k\: \R^n \ra \R^n$ defined by 
$$ g_k(x)= \tfrac 1{r_k} g(r_k x) \forr x\in \R^n$$ 
converge to the linear map  $L=Dg(0)$ locally uniformly on $\R^n$.      
Hence the balls
$$ \tfrac 1{r_k}  D'_k =\tfrac 1{r_k}  g(D_k) = g_k(\tilde D_k)$$
Hausdorff converge 
to the closed  set $D':=L(D)$ as $k\to \infty$.  It follows that $D'$ is also a closed 
ball, possibly degenerate.  
Since every linear  transformation on $\R^n$ that maps a non-degenerate ball
to a ball is conformal, the result follows. 
\hfill $\Box$ \medskip

\no {\em Proof of Theorem~\ref{T:Rig1}.}
Let $S$ and $S'$ be Schottky sets in $\Sph^n$, $n\ge 2$,
 and $f\: S\ra S'$ 
a quasisymmetric map.  Assume that $S$ has measure zero. 
We have to show that $f$ is the restriction of a M\"obius transformation. 

We use the notation of Section~\ref{s:ext}, and   let 
 $F\co \Sph^n\to\Sph^n$ be   the equivariant quasiconformal
 extension of $f$ obtained in Proposition~\ref{prop:eqvar2}. 
We will show that $F$ is a M\"obius transformation.

The set $S_\infty$ is a union of a countable number of copies of $S$ under 
M\"obius transformations.  Since $S$ has measure zero, the same is 
true for $S_\infty$, and so the set $\Sph^n\setminus S_{\infty}$
has full measure. 
By Lemma~\ref{lem:nested} 
each point in the  set  $\Sph^n\setminus S_{\infty}$ 
is contained in a sequence of closed balls $D_k$ with $\diam(D_k)\to 0$
such that each ball $D_k$ 
is an image of a ball in the collection 
 $\{ \bar B_i: i\in I\}$ under a  M\"obius transformation
in  $\Gamma_S$. Since $F$ maps peripheral spheres of $S$ to peripheral 
spheres of $S'$ and is equivariant, it follows that  $D_k':=
F(D_k)$ is a ball  for each $k\in \N$.   

Since $F$ is quasiconformal, there exits a set $N\sub \Sph^n$ of measure 
zero such that $F$ is differentiable with invertible derivative 
 at each point in $\Sph^n\setminus N$ \cite[Ch.~4]{Va}. 
 Lemma~\ref{lem:conf} implies that for each point in $\Sph^n\setminus(S_\infty\cup N)$
the map $F$ is differentiable with a derivative that  is 
an invertible conformal linear map. 
Since $\Sph^n\setminus(S_\infty\cup N)$ has full measure, the map $F$ 
is $1$-quasiconformal, and hence a M\"obius transformation. 
\hfill $\Box$ 

\medskip
Very similar arguments can be found in the proofs of Lem\-mas~3.14 and 3.15 in \cite{fri1}.

\section{Schottky sets of positive measure}
\label{s:pos} 

\no
We identify $\Sph^2$ with the extended complex plane $\oC=\C\cup\{\infty\}$, 
and denote by  $z$ a variable point in  $\oC$. By definition
a  {\em Beltrami coefficient}
is an essentially bounded complex measurable function $\mu$ on $\oC$ with  $\Vert \mu \Vert_\infty< 1$.
Each Beltrami coefficient $\mu$ defines a conformal class of 
measurable Riemannian metrics $ds^2$ on $\oC$
 (and hence a unique conformal structure) by setting
$$ d s^2= \la(z) |dz+\mu(z)d\bar z|^2, $$
where $\la$ is an arbitrary measurable function on $\oC$ that is positive almost everywhere. 
To an arbitrary  orientation preserving quasiconformal map
$F\: \oC  \ra \oC$ one can associate a Beltrami coefficient  
$\mu_F$  defined as 
$$ \mu_F= F_{\bar z}/F_z$$ 
 for almost every $z\in \oC$,  
where $F_{\bar z}=\frac{\partial F}{\partial \bar z}$ and $F_z=\frac{\partial F}{\partial z}$. If $F$ is orientation reversing, then we define $ \mu_F=   \mu_{\bar F}$, where $\bar F(z)=\overline {F(z)}$. 
If $ds^2$ is a measurable Riemannian metric associated with  a Beltrami coefficient $\mu$ and $F$ is a quasiconformal map on $\oC$, then the pull-back $F^*(ds^2)$ of $ds^2$ by $F$ is well-defined and lies in the conformal class determined by a Beltrami coefficient $\nu$,
called the {\em pull-back} of $\mu $ by $F$, written $\nu=F^*(\mu)$.  We have 
$$ F^*(\mu)= \frac {\mu_F +(\mu\circ F) \overline{F_z}/F_z  }
{ 1+  \overline{\mu_F}\, (\mu\circ F)  \overline{F_z}/F_z }
\quad \text{or}\quad  
  F^*(\mu)= \frac {\mu_F +(\overline{\mu\circ F}) 
 F_{\bar z}/\overline {F_{\bar z }}}
{ 1+  \overline{\mu_F}\, (\overline {\mu\circ F})  F_{\bar z}/
\overline {F_{\bar z} }  }
  $$ 
depending on whether $F$ is orientation preserving or reversing. 
In particular,   if $F$ is orientation reversing, then 
 $F^*(\mu)={\bar F}^*(\tilde  \mu)$, where $\tilde \mu(z)= \overline {\mu (\bar z)}$.  Note that $\mu_F$ is the pull-back by $F$ of the Beltrami coefficient $\mu_0\equiv 0$ that  defines the standard conformal structure on $\oC$. The pull-back operation on Beltrami coefficients has the usual functorial properties: If $F$ and $G$ 
 are quasiconformal maps on $\oC$ and $\mu$ is a Beltrami coefficient, then $(F\circ G)^*(\mu)= G^*(F^*(\mu))$. 
 
The Measurable Riemann Mapping Theorem says 
that  for a given Beltrami coefficient $\mu$,  there exists a quasiconformal  mapping $F$ on $\oC$ with $\mu_F=\mu$. The map $F$ is uniquely determined up to 
post-composition by a M\"obius transformation. 

Let $\Gamma$ be a group of M\"obius transformations. We say that 
a Beltrami coefficient $\mu$ is invariant 
under $\Gamma$ if  $\ga^*(\mu)=\mu$ for every $\gamma\in \Gamma$. This is equivalent to
\begin{equation}\label{Gpinvariant}
\mu =(\mu\circ\gamma)\cdot\frac{\overline{\gamma_{z}}}
{\gamma_z}\quad
\text{or} \quad
\mu =(\overline{\mu\circ\gamma})\cdot\frac {\gamma_{\bar z}}
{\overline{\gamma_{\bar z}}}
\end{equation}
almost everywhere on $\oC$ depending on whether $\gamma$ is  orientation preserving  or 
orientation reversing. 

\begin{lemma}  \label{L:gpinv}Let $\Gamma$ be a group of M\"obius transformations on 
$\oC$, and 
$F\: \oC \ra \oC$ a quasiconformal map with a Beltrami coefficient 
$\mu_F$ invariant under $\Gamma$. 
Then $F$ conjugates $\Gamma$ to a group of M\"obius transformations, 
i.e., $F\circ \gamma \circ F^{-1} $ is a M\"obius transformation 
for every $\gamma \in \Gamma$.   
\end{lemma} 

\no {\em Proof.} For every $\ga \in \Gamma$  
the map $F\circ\gamma\circ F^{-1}$ is quasiconformal. It pulls back the Beltrami coefficient $\mu_0\equiv 0$ defining the standard conformal structure to itself. This  follows from a straightforward computation using the functorial properties of the pull-back operation and the invariance of $\mu_F$ under $\Gamma$. 
This implies that   $F\circ\gamma\circ F^{-1}$ is conformal or anti-conformal 
depending on whether $\ga$ is  orientation preserving or not.  
Hence   $F\circ\gamma\circ F^{-1}$ is a M\"obius transformation
for every $\ga\in \Gamma$.  
\hfill $\Box$ \medskip

\begin{lemma}\label{L:SchtoSch}  Let $S$ be a Schottky set in $\oC$,  
 and
$F\: \oC \ra \oC$ a quasiconformal map with a Beltrami coefficient 
$\mu_F$ invariant under the Schottky group $\Gamma_S$ associated with 
$S$. Then $S'=F(S)$ is a Schottky set.   
\end{lemma}

\no {\em Proof.} Let $R$ be a reflection in one of the 
peripheral circles  $\Sigma$ of $S$. Then $R'=F\circ R \circ F^{-1}$ 
is a M\"obius transformation by the previous lemma. Since $R'$ 
is orientation reversing  and has infinitely many fixed points, it has to be a reflection in a circle $\Sigma'\sub 
\oC$. Under the map $F$ the fixed point set $\Sigma$ of $R$ 
corresponds to the fixed point set $\Sigma'$ of $R'$. 
Hence $F$ maps each peripheral circle of $S$ to a circle. 
It follows that $S'=F(S)$ is a Schottky set. 
\hfill $\Box$ \medskip 

\begin{lemma}\label{L:ident}
Suppose $U$ is an open subset in $\R^n$ with $0\in U$, and $f\: U\ra \R^n$ is a mapping that is differentiable at $0$. If there exists a set $S\sub U$ that has a Lebesgue density point at $0$ such that $f|S=\id_S$, then $Df(0)=\id_{\R^n}$.
\end{lemma}

\no {\em Proof.} For each $\eps>0$ the set 
$$M_\eps=\{ s/|s|: s\in S \text { and } 0<|s|<\eps\} \sub \Sph^{n-1}$$ 
is dense in $\Sph^{n-1}$; for otherwise,  a truncated cone with vertex at $0$ would be contained in the complement of $S$, and so $0$ would not be a Lebesgue density point of $S$. 
Hence if $\zeta\in \Sph^{n-1}$ is arbitrary, there exists a sequence $(s_k)$ in $S\setminus\{0\}$ such that $|s_k|\to 0$ and 
$s_k/|s_k|\to \zeta$ as $k\to \infty$. Setting $L=Df(0)$ and using our assumptions we obtain
$$ L(\zeta)=\lim_{k\to \infty} L(s_k/|s_k|)= 
\lim_{k\to \infty} \frac 1 {|s_k|} \big (f(s_k)+ o(|s_k|))=\zeta.$$
It follows that $L=\id_{\R^n}$ as desired. 
\hfill $\Box$ \medskip

\no {\em Proof of Theorem~\ref{T:Rig2}.} 
Let $S$ be a Schottky set in $\Sph^2$ which we identify with $\oC$. 
If $S\sub \Sph^2$ has  measure zero, then $S$  is rigid 
by Theorem~\ref{T:Rig1}. 

Conversely, suppose that $S$ has positive measure. 
Let $\nu$ be a non-trivial Beltrami coefficient supported on $S$, say $\nu\equiv1/2$ on $S$ and $\nu\equiv 0$ elsewhere.
Let $\Gamma_S$ be the Schottky group associated with $S$.  As was pointed out after the proof of Lemma~\ref{lem:eqvar1}, the sets $U(S)$, $U\in \Gamma$, form a measurable 
 partition of 
 $$ S_\infty= \bigcup_{U\in\Gamma_S}U(S).$$
  This implies that if we put   
$$
\mu(z)= \left\{ \begin{array}{cl} U^*(\nu)(z) & \mbox{if $z\in U(S)$ for some $U\in \Gamma$,}\\
0 & \mbox{otherwise,} 
\end{array} \right.   $$
then $\mu$ is an almost everywhere well-defined Beltrami coefficient invariant under $\Gamma_S$.   
Let $F\: \oC\to\oC$ be a quasiconformal map with Beltrami coefficient $\mu_F=\mu$ almost everywhere.   Then $F$ is quasisymmetric, and $F(S)$ is a Schottky set  by  Lemma~\ref{L:SchtoSch}. Moreover,  $F$ does not agree with any M\"obius transformation on $S$. For suppose it did. Then post-composing 
$F$ by a M\"obius transformation if necessary, we may assume that $F|S=\id_S$.  Then by Lemma~\ref{L:ident}, the map 
$DF(z)$ is the identity for almost every $z\in S$. This implies that $\mu_F(z)=0$ for almost every $z\in S$. This contradicts the fact that $\mu_F(z)=\nu(z)= 1/2$ for almost every  $z\in S$, because $S$ has positive measure. 
 
This shows that $S$ is not rigid. 
 \hfill $\Box$ \medskip  
   
We now give an example of a Schottky set in $\Sph^n$, $n\ge 2$,
 that has  empty
interior and is not rigid. For simplicity we work with Schottky sets in $\R^n$ (defined in the obvious way). A Schottky set in $\Sph^n$ can be obtained by adding the point at infinity.
 A similar example is contained in~\cite{KAG86} and is originally due to Apanasov~\cite{bA78}.

\begin{example}\label{ex:nonrig} 
Let $K$ be a compact set in $\R$ of positive measure, but
with no interior points. For example,  a ``thick" Cantor set will have this
property. We may assume that $0\in K$.
 The complement of $K$ can be written as
  $\R\setminus K= \bigcup_{k\in \N}I_k$, 
where the sets $I_k$ are pairwise disjoint  open intervals. 
There exists a unique absolutely continuous function $h\: \R\ra \R$ with 
$h(0)=0$ such that $h'(x)=2$ for almost every $x\in K$ and 
$h'(x)=1$ for every $x\in \bigcup_{k\in \N}I_k$.
Obviously, $h$ is a bi-Lipschitz homeomorphism of $\R$ onto itself which is different
from the identity map on $K$ and is a translation if restricted to any of the 
intervals $I_k$. Define a homeomorphism $F\:\R^n\ra \R^n$, $n\ge 2$,
by 
$$ F(x_1, \dots, x_n):= (h(x_1), x_2, \dots, x_n)$$
for $(x_1, \dots, x_n)\in \R^n$.
Then $F$ is a bi-Lipschitz homeomorphism of $\R^n$ onto itself, and is a translation if restricted to any of the slabs $M_k=I_k\times \R^{n-1}$. 
Each slab $M_k$ can be filled out with open balls such that no interior 
remains; more precisely, for each $k\in \N$ there exist pairwise disjoint
open balls $B_{kl}\sub M_k$, $l\in \N$, such that $M_k\setminus 
\bigcup_{l\in \N} B_{kl}$ has empty interior. 
Then $S=\R^n\setminus\bigcup_{k,l\in \N}B_{kl}$
is a Schottky set in $\R^n$ without interior points. Moreover, 
since $F$ restricted to the slab $M_k$ is a translation 
and each ball $B_{kl}$ lies in $M_k$, it follows that 
$B'_{kl}:=F(B_{kl})$ is a ball for all $k,l\in \N$. Hence
$$ S'=F(S) = \R^n\setminus\bigcup_{k,l\in \N}B'_{kl}$$ 
is a Schottky set.  As the restriction of a bi-Lipschitz homeomorphism
the map 
$f=F|S$ is a quasisymmetry and maps the Schottky set $S$ to the 
Schottky set $S'$. Moreover,  
$f$ is not the restriction of a M\"obius transformation. 
Indeed, suppose that $f=U|S$ for some M\"obius transformation 
$U$. By construction of $S$ we have 
$$\{0\}\times \R^{n-1}\sub K \times \R^{n-1}\sub S.$$
Since $h(0)=0$, this implies that $U$  is the  
identity on $\{0\}\times \R^{n-1}$; but $U$ has to preserve orientation
and so  $U$ is the identity map. Hence $f$ is the identity map 
on $S$ which implies that  $h$ is the identity map on 
$K$. Since this is not the case, we get a contradiction 
showing that $S$ is not rigid. 

\end{example}

\section{Rigid Schottky sets of positive measure}
\label{s:posrig} 
\no
In this section we give an example of a Schottky set in 
$\Sph^n$, $n\ge 3$,
 that  has  positive measure  and is rigid. 
 We first discuss some terminology.  In this section it is convenient  
to identify $\Sph^n$ with $\R^n\cup\{\infty\}$ (equipped with the chordal metric) 
via stereographic projection. 

Let $A$ be a subset of $\Sph^n=\R^n\cup\{\infty\}$ with $0\in A$.
We say that a set $A_\infty\sub \Sph^n$ is a {\em weak tangent}
of $A$ (at $0$), if it is closed and if there exists a sequence $(r_k)$ 
of positive numbers tending to $0$ such that $A_k\to A_\infty$, 
where 
$$ A_k=\tfrac 1{r_k}A= \{\tfrac 1{r_k}x: x\in A\}. $$
Here  we use the convention that $\la\cdot \infty =\infty $ for all $\la>0$. 
So a weak tangent of $A$ is a closed set that we obtain by ``blowing up" $A$ at the origin in a suitable sense.  Every set $A$ with $0\in A$ has a weak tangent, because for every sequence $(r_k)$ of positive numbers with $r_k\to 0$, the sequence of sets $A_k=\frac 1 {r_k} A$ subconverges. 
Every  weak tangent of $A$ contains the point $\infty$ unless $A=\{0\}$.

Our notion of a weak tangent is 
suitable for our purposes and is a  variant  of  similar concepts  in the literature.

   Let $D$ be a region in $\Sph^n$ (i.e., an open and connected subset of $\Sph^n$),  and let $T$ be a subset of $D$ whose complement in $D$ is a union of at least three disjoint open balls such that the  closure of each ball is  contained in $D$. Such a set $T$ will be called a \emph{relative Schottky set in} $D$. 
The boundaries of the balls in the complement of $T$ in $D$ are referred to as peripheral spheres.

If $\Sigma$  is a peripheral sphere of $T$, then 
$T\setminus \Sigma$ is path-connected. Indeed, to connect two points $x,y\in T\setminus \Sigma$, one first takes an arc in $D\setminus \Sigma$ 
that consists of finitely many spherical  geodesic segments and  joins  
$x$ and $y$.  Then one proceeds similarly as 
in  the proof of Lemma~\ref{lem:ballconn} to ``correct" 
$\gamma$  on suitable subarcs to create a path 
$\tilde \gamma $ in $T\setminus \Sigma$ joining 
$x$ and $y$. 

Moreover, if $\Sigma$ is any topological 
$(n-1)$-sphere in $T$ such that $T\setminus \Sigma$ is connected, then $\Sigma$ is a peripheral sphere 
of $T$. This follows from the second part of the proof 
of Proposition~\ref{prop:peri} applied to $S=T$. (Note 
that  $K'$ has to meet one of the complementary 
components $B$ of $T$ in $D$, for otherwise the non-empty set $\partial K'$ would be contained in 
$(\Sph^n\setminus D)\cap \Sigma=\emptyset$.)   

This shows that if  $\Sigma$ is any topological 
$(n-1)$-sphere in $T$, then  $T\setminus \Sigma$ is connected if and only if $\Sigma$ is a peripheral sphere of $T$. In particular, every homeomorphism 
 between relative Schottky sets has to take peripheral spheres to peripheral spheres.
  
A relative Schottky set $T$ in $D$ is called \emph{locally  porous at} $x\in T$ if  there exist an open neighborhood $U$ of $x$ and constants $C\ge 1$ and $\rho_0>0$ with the property that for each $y\in T\cap U$ and each $r$ with $0<r<\rho_0$ there exists a complementary component $B$ of $T$ in $D$  with    $B(y,r)\cap B \ne \emptyset$ and   $r/C\leq\diam(B)\leq Cr$.
 If this is true for each point $x\in T$, then we call $T$  \emph{locally  porous}. 
 A locally porous relative Schottky set cannot have Lebesgue density points, and hence is a set of measure naught. 


Our first goal is the proof of the following theorem.
\begin{theorem}\label{T:relSch}
Let $n\in\N,\ n\geq 3$, $T$ and $T'$ be relative
 Schottky sets in regions $D,D'\sub \Sph^n$, respectively, and $\psi\: T\to T'$ a  quasisymmetric map. If $T$ is locally  porous, then  $\psi$ is the restriction of a M\"obius transformation to $T$.
\end{theorem}

 We need the following lemmas. 

\begin{lemma}\label{L:invpor}
Let $T$ be a relative Schottky set in a region 
$D\sub \Sph^n$, $n\in \N$, and $F\: \Sph^n\ra \Sph^n$ a M\"obius transformation. Then  $T'=F(T)$ is a relative  Schottky set in 
$D'=F(D)$. 
If in addition $T$ is locally  porous at $x\in T$, then $T'$ 
 is locally  porous at $x'=T(x)$. 
\end{lemma}

\no{\em Proof.}  It is clear that $T'$ is a relative Schottky set in $D'$. 
Assume $T$ is locally  porous at $x\in T$. We have to show 
that $T'$ is  locally  porous at $x'=T(x)$. 
To see this we use the following general fact whose proof is left to the reader: Suppose  $G\: \Sph^n \ra \Sph^n$ 
is an $\eta$-quasisymmetric map, and $M,N\sub \Sph^n$ are two sets with 
$M\cap N\ne \emptyset$ and 
$$ (1/C)\diam(M)\le\diam(N)\le   C\diam(M), $$
where $C\ge 1$. Then for $M'=G(M)$ and $N'=G(N)$, we have 
$$ (1/C')\diam(M')\le\diam(N')\le   C'\diam(M'), $$
where $C'\ge 1$ only depends on $C$ and $\eta$. 
In other words, if $G$ is a quasisymmetric map, then the images under $G$  of two intersecting sets that  have comparable size will also have comparable size, quantitatively. 

The claim now follows if we apply this statement to $G=F$ and to the sets 
$M=B(y,r)$ and $N=B$ appearing  in the definition 
of local  porosity.  We leave the details to the reader. 
\hfill$\Box$\medskip

\begin{lemma}\label{L:conSch} Suppose  $T_\infty$ is a weak tangent  of a relative Schottky set $T$ in a region $\Omega\sub \Sph^n$ with $0\in T$. Then the complementary components of 
$T_\infty$ in $ \Sph^n$ are open balls. \end{lemma}

In particular, $T_\infty$ is a 
 Schottky set if it has  at least three such components.   

\medskip
\no{\em Proof.}  We can write $T=\Omega\setminus \bigcup_{i\in I} D_i$, where 
the sets $D_i$, $i\in I$, form  a family of disjoint open balls in  
$\Omega$. 
There exists a sequence $(r_k)$ of positive numbers tending to $0$ such that 
$T_k\to T_\infty$, where $T_k=\frac 1{r_k} T$. 

Now  let $x\in \Sph^n\setminus T_\infty$ be arbitrary. Then $x\in \Sph^n\setminus T_k$ for large $k$. (Note that here  we use that $T_\infty$ is a closed set). 
Since $\infty\in T_\infty$, we have $x\ne \infty$. Moreover,   since $0\in \Omega$, and so $x\in
\frac 1{r_k} \Omega$ if $k$ is large, we can find $k_0\in \N$ 
and $i_k\in I$ for $k\ge k_0$ such that $x\in B_k:=  \frac 1{r_k}
D_{i_k}$. The sequence of balls $(B_k)$ subconverges to a closed ball $B$; keeping the same notation for this subsequence
for convenience, we may assume $B_k\to B$. 
Then $x\in B$. 
  Suppose $y\in \inte(B)$. By the last part of 
Lemma~\ref{L:ballconv}, there exists $\delta>0$ such that 
$B(y,\delta)\sub B_k$ and so $\dist(y,T_k)\ge \delta$ for large $k$. Hence $y\in \Sph^n\setminus T_\infty$. This shows that $\inte(B)\sub  
\Sph^n\setminus T_\infty$.   By Lemma~\ref{L:ballconv} we also have $\partial B_k\to \partial B$. Since $\partial B_k \sub T_k$, it follows that $\partial B\sub 
T_\infty$. We conclude  that  the open ball $\inte(B)$ is the connected component 
of the complement of $T_\infty$ containing $x$.
 Since $x\in \Sph^n\setminus T_\infty$ was arbitrary, the claim follows. 
\hfill $\Box$ \medskip

\begin{lemma}\label{L:porSch} Suppose $T$ is a relative Schottky set 
that is locally  porous at $0\in T$, and  $T_\infty$ 
is a  weak tangent  of $T$. Then $T_\infty$ is a Schottky set 
that is  locally  porous at every point $x\in T_\infty\setminus\{\infty\}$. 

 In particular, $T_\infty$  has measure zero. \end{lemma}

Actually, one can show that $T_\infty$ is also locally  porous at $\infty$,
but we do not need this fact for the desired conclusion that $T_\infty$ has 
measure zero.  

\medskip 
\no{\em Proof.} We use notation as in Lemma~\ref{L:conSch} and its proof. 
Near each point in $\Sph^n\setminus \{\infty\} =\R^n$ the Euclidean 
metric and the chordal metric are bi-Lipschitz equivalent. 
Therefore, we can use our assumption that $T$ is locally  porous 
at $0$ and  derive the desired conclusion that $T_\infty$ is    locally  porous at every point $x\in T_\infty\setminus\{\infty\}$ by 
using the Euclidean metric instead of the chordal metric.  For the rest  of the proof all metric notions refer to the Euclidean metric on $\R^n$.

 The neighborhood and the constants in the definition of local porosity of $T$ at $0$  will be denoted by $U$, $C$ and $\rho_0$, respectively.  Let $x$ be an arbitrary point in $T_\infty\setminus\{\infty\}$ and $R>0$. The point $x$ is the limit of a sequence $(x_k)$ such that $x_k\in T_k$. 
For sufficiently large $k$ we have  $r_kR< \rho_0$ ,  
and $r_k x_k\in T\cap U$. 
Using the local  porosity of $T$ it follows that then 
there  exists $i_k\in I$  such that  $$ B(r_k x_k,  r_k R/2)\cap D_{i_k}\ne \emptyset\quad
\text{and}\quad Rr_k/C'\leq\diam (D_{i_k})\leq C'R r_k, 
$$
where $C'=2C$. If we define $B_k:=\frac 1{r_k} D_{i_k}$, then $B_k$ is a complementary component of $T_k$, and the previous statements  translate to 
\begin{equation}\label{E:sph}
 B( x_k,  R/2)\cap  B_k\ne \emptyset
\quad\text{and}\quad  
R/C'\leq\diam ( B_k)\leq C'R. 
\end{equation}

Passing to an appropriate subsequence if necessary, we may assume that $B_k\to B'$, where $B'$ is a closed ball.  Let $B=\inte(B')$. Then by 
\eqref{E:sph} we have 
\begin{equation}\label{E:sph2}
 B(x, R)\cap B\ne \emptyset \quad\text{and} \quad
R/C'\leq\diam (B)\leq C'R.
 \end{equation}
  Moreover, the argument in the proof of Lemma~\ref{L:conSch} shows that $B$ is a complementary component of 
 $T_\infty$ in $\Sph^n$. 
 Since $x\in T_\infty\setminus\{\infty\}$ and $R>0$ in  \eqref{E:sph2} are 
 arbitrary, the local  porosity of $T_\infty$ at each point different
  from $\infty$ now follows.  Moreover,  \eqref{E:sph2} also shows that $T_\infty$ has infinitely many 
complementary components and is hence a Schottky set by 
Lemma~\ref{L:conSch}. 

Finally,  $T_\infty$ is a set of measure zero, because  
$T_\infty$ cannot have any Lebesgue density points except possibly the point
$\infty$.
\hfill$\Box$\medskip

\no{\em Proof of Theorem~\ref{T:relSch}.}
Let $\Sigma$ be any peripheral sphere of $T$ and $\Sigma'=\psi(\Sigma)$ be the corresponding peripheral sphere of $T'$. The restriction $\phi=\psi|\Sigma$ of $\psi$ to  $\Sigma$ is a quasiconformal map between $(n-1)$-dimensional spheres. Therefore, at almost every point of $\Sigma$ (with respect to  spherical $(n-1)$-dimensional measure) the map $\phi$ is differentiable with invertible derivative.

We want  to show that at every such point $x_0\in \Sigma$ the derivative is a conformal map (i.e., a scalar multiple of an isometry). By composing with M\"obius transformations we can assume that
$\Sigma=\Sigma'=\R^{n-1}\cup\{\infty\}$, 
$\psi(x_0)=\phi(x_0)=x_0=0$ and $\psi(\infty)=\phi(\infty)=\infty$. 
Here we make the identification $\Sph^n=\R^n\cup\{\infty\}$, 
and consider $\R^{n-1}$  as a subset of $\R^n$
in the usual way. Note that by Lemma~\ref{L:invpor} our assumptions on $T$  
are  not affected by such auxiliary M\"obius transformations.

We  extend the  quasisymmetric map $\psi\: T\ra T'$ 
to a quasiconformal map $F\: D\ra D'$. 
The existence 
of such an extension follows from the same method as in the proof of
Proposition~\ref{prop:ext}.     

There exists a sequence $(r_k)$ of positive numbers tending to $0$ such that 
$$ T_k=\tfrac1{r_k}T\to T_\infty \quad\text{and}\quad T_k':= \tfrac1{r_k}T'\to T'_\infty,$$ where $T_\infty$ and $T'_\infty$ are weak tangents of $T$ and $T'$ (at $0$), respectively. 

 Consider the maps $F_k$ defined by 
$F_k(x)=F(r_kx)/r_k$ for $k\in \N$. Since $0\in D$ (after applying the auxiliary 
M\"obius transformation discussed above), the maps $F_k$ are eventually
defined on every ball $B(0,R)$, $R>0$, (with respect to the Euclidean metric on $\R^n$) and map  $B(0,R)$ into $\R^n$ .  Moreover,       the sequence of maps $(F_k)$ is uniformly quasiconformal, i.e., there exists $H\geq 1$ such that each map $F_k$ is $H$-quasiconformal. Also, $F_k(0)=0$,  and if 
$e_1=(1,0\dots,0)\in \R^n$, 
$$
\lim_{k\to\infty}  F_k(e_1)=\lim_{k\to\infty}  \phi(r_k e_1)/r_k=D\phi(0)(e_1)\neq 0. 
$$
Using standard compactness arguments for quasiconformal maps (see \cite[Sect.~21]{Va}),  we conclude that there exists a subsequence of $(F_k)$ that  converges locally uniformly to a quasiconformal map $F_\infty$
on $\R^n$. For convenience of notation we continue indexing this subsequence by $k$. By applying a similar argument to the 
inverse maps $G_k:=F_k^{-1}$, we may in addition assume 
that $G_k\to G_\infty$ locally  uniformly on $\R^n$, where $G_\infty$ is a quasiconformal map on $\R^n$. By putting $F_\infty(\infty)=G_\infty(\infty)=\infty$ we can extend these maps to quasiconformal maps on $\Sph^n$.   
Then  $F_\infty\circ G_\infty=G_\infty\circ F_\infty=\id_{\Sph^n}$, and so 
$F_\infty$ and $G_\infty$ are inverse maps of each other. 
 
We claim that $F_\infty(T_\infty)=T'_\infty$. Since $G_\infty=F_\infty^{-1}$, this is equivalent to the inclusions  $F_\infty(T_\infty)\sub T'_\infty$ and
 $ G_\infty(T_\infty')\sub T_\infty$.
By symmetry it is enough to show the first inclusion. So let $x\in T_\infty$ be arbitrary. If $x=\infty$, then $F_\infty(x)=\infty\in T'_\infty$. If $x\in 
T_\infty\setminus\{\infty\}$,  then there exists a sequence $(x_k)$ with $x_k\in T_k$ for $k\in \N$ and $x_k\to x\in \R^n$. Since $F_k\to F_\infty$ locally uniformly on $\R^n$, it follows that $F_k(x_k)\to 
F_\infty(x)$. 
On the other hand, $F_k(x_k)\in \frac1 {r_k} F(T)= T'_k$.
Hence $F_\infty(x)\in T'_\infty$. 

According to Lemma~\ref{L:porSch}, the set
$T_\infty$ is a  Schottky set of measure  zero. Using  Lemma~\ref{L:conSch},
and the fact that $T'_\infty=F_\infty(T_\infty)$, we see that   $T'_\infty$
is also a Schottky set.    The map $F_\infty$ is quasiconformal, and hence   quasisymmetric, since $n\ge 2$.  It follows that we can apply Theorem~\ref{T:Rig1} to conclude that $F_\infty$ agrees with a 
M\"obius transformation on $\R^{n-1} \sub T_\infty$. 
Since $F_\infty(\infty)=\infty$, $F_\infty (0)=0$, and $F_\infty(\R^{n-1})\sub 
\R^{n-1}$, the map  $F_\infty|\R^{n-1}$ has to be  a conformal linear map. 
 
On the other hand, it follows from the definitions of $F$ and  $F_\infty$ that  
$F_\infty|\R^{n-1}=D\phi(0)$ which proves the desired statement that $D\phi(x_0)=D\phi(0)$ is conformal. 

Since this is true for almost every point $x_0\in \Sigma$, the map  
$\phi$ is a $1$-quasiconformal map between the peripheral spheres $\Sigma$ and $\Sigma'$.  Hence it is the restriction of a M\"obius transformation on  $\Sph^n$ to $\Sigma$. 

 If $B$ is a complementary component of $T$ in $D$, then $\Sigma=\partial B$ is a peripheral sphere of $T$. Since $\Sigma'=\psi(\Sigma)$ is a peripheral sphere of $T'$, there exists a  corresponding complementary component $B'$ of $T'$ such that $\partial B'=\Sigma'$. By what we have seen in  the first part of the proof  there exists a M\"obius transformation that agrees with $\psi$ on $\Sigma$ and  maps $\bar B$ to $\bar B'$. Using such M\"obius transformations we can extend  the original map $\psi$ to each  complementary  component of $T$ in $D$, to obtain a bijection  $\Psi\: D\ra D'$. As in the proof of Proposition~\ref{prop:ext} one can show that this extension is a quasiconformal homeomorphism.
On each of the complementary components of $T$ in $D$ the map 
$\Psi$ agrees with a M\"obius transformation. Moreover, since $T$ is locally  porous,  it has measure zero, and so the complentary components of $T$ in $D$
form a set of full measure in $D$. It follows that $\Psi$ is $1$-quasiconformal.   Since $n\ge 3$, we can apply Liouville's Theorem, and so the map $\Psi$ and hence also $\psi$ is the restriction of a 
M\"obius transformation.      
\hfill$\Box$\medskip

After these preparations we are ready to construct 
rigid Schottky sets with positive measure in dimension $n\ge 3$.  
\medskip 

\no{\em Proof of Theorem~\ref{T:exposrig}.}
Let $D$ be a region in $\Sph^n$ (i.e., an open and connected set) that  is dense and whose boundary has positive measure. For example, one can take the complement of a  ``thick" Cantor set for $D$. We want to show that $D$ contains a locally  porous relative Schottky set $T$. The existence of such a set implies the statement by Theorem~\ref{T:relSch}.  Indeed, the set $S=T\cup\partial D$ is a Schottky set of positive measure. Every quasisymmetric map $f$ of $S$ onto any other Schottky set restricts to $T$ as a quasisymmetric map onto another relative Schottky set, and is therefore the  restriction of a M\"obius transformation to $D$. Since $D$ is dense in $\Sph^n$, we conclude that $f$ is the restriction of a M\"obius transformation  to $S$.

In order to construct a locally  porous relative Schottky set in $D$, we proceed as follows.  Consider the subset $N_1$ of $D$ defined by 
$$N_1=\{x\in D\: {\rm dist}(x,\partial D)\geq1\}.$$  Let $A_1$ denote a  maximal $1$-separated subset of $N_1$, and let $D_1$ be the set obtained from $D$ by removing the union of all disjoint open balls with radii $1/4$ centered  at  elements of $A_1$ . 
Inductively, if $k\in \N,\ k\geq 1$, let 
$$N_{k+1}=\{x\in D_k\: {\rm dist}(x,\partial D_k)\geq1/2^{k}\},$$ let $A_{k+1}$ be a maximal $1/2^k$-separated subset of $D_k$, and let $D_{k+1}$ be the set obtained from $D_k$ by removing the union of all disjoint open balls with radii $1/2^{k+2}$ centered at elements of
$A_{k+1}$ . The sets $D_k$ form a monotonically decreasing  sequence of subsets of $D$, and their intersection is by construction a relative Schottky set $T$ in $D$. 

To show that $T$ is locally  porous, 
let $x\in T$ be arbitrary.  Define  $d=\dist(x, \partial D)$, $U=B(x,d/2)$, and $\rho_0=d/4$, 
 and suppose  $y\in T\cap U$ and  $r$ with $0<r<\rho_0$ are arbitrary.  

 By construction of $T$ there exist infinitely many complementary components of $T$ in $D$ intersecting $B(y,r)$.  Among all such components,      we can choose one, say $B_0$, with largest diameter. 
 Since $$\dist(B(y,r), \partial D)\ge d/4\ge r,$$ the   construction of $T$ shows that  $\diam (B_0)\geq 
 c_1r$, where $c_1>0$ is an absolute constant. In general an inequality of this type  will not be true in the other direction, because $\diam( B_0)$ can be much larger than $r$. 
To obtain a complementary component that intersects $B(y,r)$ and has diameter comparable to $r$, we take  the  second largest 
complementary component that meets $B(y,r)$.  
 
 More precisely, 
let $B$ be a complementary component of $T$ in $D$ different from $B_0$ with  $B(y,r)\cap  B \ne 0$  that  has largest diameter among all such components. 
   Since $\dist ( B,  B_0)\le 2r$, by construction of
   $T$ we have $\diam( B)\le c_2r$ 
   where $c_2>0$ is an absolute constant. 
  On the other hand, it is clear that  $\diam( B)\ge c_3 r$ for some absolute constant $c_3>0$.  
This implies  that the complemetary component $ B$ is the desired one in the condition for the local  porosity where  
  $C$   can be taken as an absolute constant. So $T$ is locally  porous.  The proof is complete. 
\hfill$\Box$\medskip

The construction in the previous proof should be compared with 
Example~\ref{ex:nonrig}.  Essentially, in this example we constructed a relative 
Schottky set in the open set $D'=(\R\setminus K)\times \R^{n-1}$. In contrast to the set $D$ in the above proof, this set is not  connected, but consists of infinitely many open slabs.  As we saw, in this case  a quasisymmetric map can agree  
 with a M\"obius transformation  on each component of $D'$ 
  without being a M\"obius 
 transformation globally.

\section{Rigidity for convex subsets of hyperbolic space}
\label{s:con} 

\no 
Fix $n\in \N$, $n\ge 2$. We denote by $\Ca_n$ the 
class of all closed convex subsets $K$ of hyperbolic $n$-space $\H^n$ with non-empty interior and non-empty totally 
geodesic boundary. So  $\partial K\sub \H^n$ consists of  a union of pairwise disjoint hyperplanes. To rule out some trivial cases, we make the additional 
assumption that there are at least three such hyperplanes in $\partial K$. 

Usually, we think of $\H^n$ in the conformal unit ball 
model. Then the boundary at infinity $\geo \H^n$
can be identified with the unit sphere $\Sph^{n-1}$. 
If $K\in \Ca_n$, then the boundary at infinity 
$\geo K\sub \Sph^{n-1}$ is a Schottky set. 
Conversely, if $S\sub \Sph^{n-1}$ is a Schottky set, 
then its hyperbolic convex hull $K\sub \H^n$ belongs 
to the class $\Ca_n$. 

Suppose  $(X,d_X)$ and $(Y,d_Y)$ are  metric spaces. A map
$f\:X\rightarrow Y$ is called a {\em quasi-isometry}  
of $X$ into $Y$ if there exist constants $\lambda\ge 1$ and $k\ge 0$ 
such that 
\begin{equation*} \label{eqn-quasi-isometric}
\frac 1\lambda  d_X(x,x')  -  k  \le d_Y( f(x),f(x')) 
\le \lambda d_X(x,x')  +  k 
\end{equation*}
for all $x,x'\in X$,  
and if for each $y\in Y$ there exists $x\in X$ such that 
\begin{equation*} \label{surj} 
d_Y(f(x),y)\le k. 
\end{equation*}
Two maps $f,g\: X\ra Y$ are said to have  {\em finite distance} if 
$$ \dist(f,g):= \sup_{x\in X} d_Y(f(x), g(x))< \infty.$$

We call a set $K\in \Ca_n$ {\em rigid} if 
for every quasi-isometry $f\: K \ra K'$ to another  set $K'\in \Ca_n$ there exists an isometry $g$ of $\H^n$ such that $f$ and $g|K$ have finite  distance. 
Note that in this case we have $g(K)=K'$
(see the proof of Theorem~\ref{T:Rig3} below).

The following proposition records some basic properties  of quasi-isometries between Gromov hyperbolic spaces and their 
induced maps on the boundary. For the definition of a Gromov hyperbolic space and its boundary see \cite{gh},  and \cite{BS} for related considerations.  
We will use mostly  notation and terminology as in \cite{BS}. 

A proper geodesic Gromov hyperbolic  space is called {\em visual}
 if there exists a basepoint $p\in X$
and a constant $k\ge0$ such that  for 
every point $x\in X$ there exists a geodesic ray $\ga$ in $X$ with initial point 
$p$ such that  $\dist(x, \ga)\le k$ (note that the definition in \cite[p.~279]{BS} 
is equivalent in this context). 

A metric $\rho$ on the boundary $\partial_\infty X$ of a Gromov hyperbolic space $X$ is called {\em visual} if there exist a point $p\in X$ and constants
 $C\ge 1$ and $\eps>0$ such that
 $$ (1/C) e^{-\eps(u\cdot v)_p} \le \rho(u,v)\le C
 e^{-\eps(u\cdot v)_p} $$
 for all $u,v\in \partial_\infty X$. Here
 $(u \cdot v)_p$ is the ``Gromov product" 
of the points $u$ and $v$ with respect to $p$ 
(\cite[p.~273]{BS}).  Visual metrics always exist on $\partial_\infty X$, and any two of them are quasisymmetrically equivalent by the identity map.

\begin{proposition}\label{prop:qisom}
Let $X$ and $Y$ be proper geodesic metric spaces 
that are Gromov hyperbolic. Then every quasi-isometry $f\:X \ra Y$ induces a quasisymmetric map 
$\tilde f\: \geo X \ra \geo Y$.  

Suppose in addition that $X$ is visual and that $\geo X$ 
is connected and contains more than one point. Then two quasi-isometries
$f,g\: X \ra Y$ have finite distance if and only if the induced maps $\tilde f, \tilde g\: \geo X \ra \geo Y$ are identical. 
\end{proposition} 

Here we think of $\geo X $ and $\geo Y$ as being equipped with fixed visual metrics. 

\medskip
\no {\em Proof.}  These statements are essentially well-known. For the first part see \cite{BS}, Section 6, in particular  Theorem 6.5. Note that the terminology in \cite{BS} is slightly different from the one employed here. 
It follows from the definitions (see \cite[Prop.~6.3]{BS}) that two quasi-isometries 
induce the same boundary maps if they have finite distance.

Now assume in addition that $X$ is visual and $\partial_\infty X$ is connected and contains more than one point.  Suppose $ f\: X \ra Y$ is a quasi-isometry.
Fix  basepoints $p\in X$ and $q\in Y$,  and use the notation 
$z'=f(z)$ for $z\in X$ and $w'=\tilde f(w)$ for $w\in \geo X$. We will show that for every $x\in X$ the location of $f(x)$ is uniquely determined up to 
uniformly bounded distance by the data $x$ and 
$\tilde f$. This will show that if $g\: X \ra Y$ is a  quasi-isometry with $\tilde g=\tilde f$, then $f$ and $g$ have finite distance. 

Now let $x\in X$ be arbitrary. In the following, $C_1, C_2, \dots$ are  constants independent of 
$x$.  Since $X$ is visual, there exists a geodesic ray in $X$, denoted $[p,u]$,  that starts at $p$ and ``ends in" (i.e., is asymptotic to) a point   
 $u\in \geo X$ 
such that 
$$ \dist (x, [p,u]) \le C_1.$$ 
Since $\geo X$ is connected and contains more than one point, there exists a point 
$v\in \geo X$ such that 
$$ \big | (u \cdot v)_p -\dist(p,x) \big| \le C_2. $$
This inequality essentially says that the rays $[p,u]$ and $[p,v]$ start to diverge near $x$, and so 
$x$ is a   ``rough  center" of the geodesic triangle $\Delta=[p,u]\cup [p,v] \cup [u,v]$.
More precisely, 
$$ \max\{ \dist (x, [p,u]), \dist (x, [p,v]), \dist (x, [u,v])\}\le C_3, 
$$
where $[u,v]$ is the geodesic line in $X$ whose ends are asymptotic to $u$
 and $v$, respectively.    
 
Let 
$\Delta'= [p', u']\cup [p', v']\cup [u',v']$. 
By geodesic stability of Gromov hyperbolic spaces (see \cite[p.~273]{BS}), 
the image $f(\Delta)$ is within bounded Hausdorff distance of the geodesic triangle $\Delta'$.  More precisely, for the Hausdorff distance $\dist_{H}$ of these sets we have $$ \dist_{H}(f(\Delta), \Delta')\le C_4. $$ 
Note that 
$$ \dist_H([p', w], [q,w])\le C_5$$ for all $w\in \geo 
Y$, where $C_5$ is independent of $w$. 
Hence for  
$\bar \Delta  =[q,u']\cup [q,v'] \cup [u',v']$ we have 
$$ \dist_{H}(\Delta', \bar\Delta)\le C_6,$$ 
and so 
$$  \dist_{H}(f(\Delta), \bar\Delta)\le C_7. $$
It follows that $x'=f(x)$ is a ``rough center" of 
$\bar\Delta$, that is, 
$$ \max\{ \dist (x', [q,u']), \dist (x', [q,v']), \dist (x', [u',v'])\}\le C_8. 
$$ 
Since rough centers of geodesic triangles in Gromov hyperbolic spaces are essentially unique, this implies that up to contolled bounded distance, the location of $x'$ is determined by $u'$ and $v'$, i.e., by the data $x$ 
and $\tilde f$ as claimed. 
\hfill $\Box$ \medskip

If $K\in \Ca_n$, $n\ge 3$, then $K$ satisfies the assumptions on $X$ as in Proposition~\ref{prop:qisom}. First, $K$  is proper and geodesic. Moreover, $K$ 
is Gromov hyperbolic as a subset of the space $\H^n$  that has this property. 
The Gromov boundary of $K$ can be identified with the boundary at infinity 
$\geo K \sub \geo \H^n \sub \Sph^{n-1}$ of $K$ in the unit ball model of $\H^n$. Since $\geo K$ is a Schottky set, this set contains more than a point and  is connected by Lemma~\ref{lem:ballconn} if $n\ge 3$. Finally,  a set  $K\in \Ca_n$, $n\ge 2$, is visual. 
To see this fix a basepoint $p\in K$. We have to show that every point $x\in K$ lies within uniform 
distance of a ray  starting at $p$ and ending  in 
$ \geo K$.
The point $x$ lies on a ray   $[p,w]$ with $w\in \geo
\H^n$. If $w\in \geo K$ we are done. Otherwise, 
$[p,w]$ meets  one of the hyperplanes $H$ isometric to $\H^{n-1}$ comprising the set $\partial K\sub\H^n$. 
Let $[u,v]$ with  $u,v\in \geo \H^n$ be geodesic line in 
$H$ with $[u,v]\cap [p,w]\ne \emptyset$.  Since $[u,v]$ is 
contained in $H\sub \partial K$, we have $u,v\in \geo H\sub \geo K$. 
Let $H'$ be a totally geodesic subspace of $\H^n$
isometric to $\H^2$ and containing $[p,w]$ and 
$[u,v]$.  By construction the point $x$ is contained 
in the interior of the (ideal) geodesic triangle in $H'$ with   
sides $[p,u]$, $[p,v]$, $[u,v]$. By thinness of geodesic triangles in $\H^2$, this means that $x$ has uniformly bounded distance to one of the geodesic rays $[p,u]$ and $[p,v]$. Since $u,v\in \geo K$,    the set  $K$ is indeed visual. 

\medskip 
\no 
{\em Proof of Theorem~\ref{T:Rig3}.}
Let $n\ge 3$, and $K\in \Ca_n$ such that  $\geo K$ is a set of measure zero.

Suppose $f\: K\ra K'$ is a quasi-isometry to a set 
$K'\in \Ca_n$. By the discussion preceeding the proof, we can apply  Proposition~\ref{prop:qisom}, and so $f$ induces a quasisymmetric boundary map 
$\tilde f: \geo K \ra \geo K'$. Since $\geo K'$ and $\geo K$ are Schottky sets, Theorem~\ref{T:Rig1}
implies that there exists a M\"obius transformation $\tilde g$ such that $\tilde g | \geo K=\tilde f$. The map $\tilde g$ is the boundary map of an
isometry $g
\:  \H^n \ra \H^n$. Hence by 
Proposition~\ref{prop:qisom}, the maps $f$ and $g|K$ have finite distance.

We have $g(K)=K'$. Indeed,  
both sets $g(K)$ and $K'$ are in $\Ca_n$; so each set is  equal to the convex hull of its  boundary at infinity; but these boundaries are equal, since  
$$\quad\quad\quad\quad\quad\quad \geo(g(K))= \tilde g(\geo K)=\tilde f(\geo K)=\geo K'.\text{$
\quad\quad\quad\quad\quad$ $\Box$} $$  

In dimension $n=2$ and $n=3$, the rigidity of 
sets in $\Ca_n$ can be completely characterized. 

\begin{theorem}\label{T:Rig4}
No set $K\in \Ca_2$ is rigid. 
\end{theorem}

As a preparation for the proof, we first discuss some standard facts about quadrilaterals.
A {\em quadrilateral} $Q=Q(z_1,z_2,z_3,z_4)$ is a closed Jordan region in $\C$ with four distinguished points $z_1, \dots, z_4$ on its boundary. It is assumed  that the order of the points $z_k$ on the Jordan curve $\partial Q$ corresponds to positive orientation. A {\em (quasi-)conformal map} $f\:
Q(z_1,z_2,z_3,z_4)\ra Q'(z'_1,z'_2,z'_3,z'_4)$ between two quadrilaterals 
is a homeomorphism between the closed Jordan regions $Q$ and $Q'$ that is (quasi-)conformal on the interior of $Q$ and has the property that 
$f(z_k)=z'_k$ for $k=1, \dots, 4$. Every quadrilateral is conformally equivalent to a unique rectangle $R=[0,M]\times[0,1]\sub \R^2\cong \C$,
where $0$, $M$, $M+i$, $i$ are the distinguished points of $R$. The number $M>0$ is called the  {\em modulus} of $Q$, denoted by $\Mod(Q)$. 
Two quadrilaterals $Q$ and $Q'$ are conformally equivalent if and only if
$\Mod(Q)=\Mod(Q')$. In general a quasiconformal map  will distort the modulus of a quadrilateral $Q$. This distortion only depends on the Beltrami coefficient of the quasiconformal map (considered as a measurable function on $\inte(Q)$). Indeed, we have the following lemma.

\begin{lemma}\label{L:quad} Suppose $Q$ is a quadrilateral, and 
$f$ and $g$ are quasi\-conformal maps on $Q$ such that $\mu_f=\mu_g$
almost everywhere on $\inte(Q)$. Then $\Mod(f(Q))=\Mod(g(Q))$.  
\end{lemma}

\no {\em Proof:} Note that $g\circ f^{-1}$ is a quasiconformal map with a Beltrami coefficient
that vanishes almost everywhere. Hence this map is a conformal map between the quadrilaterals $f(Q)$ and $g(Q)$.
\hfill $\Box$ \medskip

For every quadrilateral $Q$ one can find quasiconformal maps that distort its modulus in a non-trivial way. Indeed, let  $f$ be a conformal map of $Q$ to a rectangle $R=[0,M]\times [0,1]$, and let  $R'=[0,M']\times [0,1]$ be any other rectangle with $M'>0$. There is a unique  affine map $A$ that
takes the quadrilateral $R$ to $R'$. Then $g=A\circ f$ is a quasiconformal map between the quadrilaterals $Q$ and $R'$. In particular, if $M\ne M'$, then $\Mod(Q)\ne \Mod(g(Q))$. 

\medskip 
\no {\em Proof of Theorem~\ref{T:Rig4}: } If $\geo K$ does not contain at least four distinct points, then $\geo K$ consists of three distinct points and $K$ is an ideal geodesic triangle. Then $K$ has bounded Hausdorff distance to a ``tripod" $T\sub K$, i.e., $T$ is  a union of three  distinct geodesic rays emanating from the same point in $K$.
In particular there is a map $g\: K \ra T$ that is the identity on $T$ and moves every point by only a bounded amount. 

 Obviously, there are quasi-isometries $f\: T\to T$ that do not 
 have finite distance to any isometry of $\H^2$ restricted to $T$; for example, such maps can be  obtained by stretching the legs 
 of the tripod by a factor $\la \ne 1$.  Then $f\circ g$ is a quasi-isometry from $K$ to $K$ that does not have finite distance 
 to any isometry on $\H^2$ restricted to $K$. Hence 
 $K$ is not rigid. 
 
For the remaining case we can assume that $\geo K\sub \partial \D$ contains four distinct points, say $z_1, z_2, z_3, z_4$, where the numbering is such that the points follow each other in positive orientation on $\partial \D$.
Here we identify $\H^2$ with the open unit disc $\D$ in $\C$ equipped with the hyperbolic metric. Then $Q=K\cup \geo K\sub \bar \D$  with the distinguished points $z_1, z_2, z_3, z_4$ is a quadrilateral. 
Fix a  Beltrami coefficient  $\nu$ on $\inte(Q)$
so that every quasiconformal map $h$ on $Q$ with $\mu_h=\nu$ almost everywhere  on $\inte(Q)$
distorts the modulus of $Q$ in a non-trivial way, i.e., $\Mod(Q)\ne \Mod(h(Q))$. This is possible by Lemma~\ref{L:quad} and the discussion following this lemma. We will use 
$\nu$ to obtain a non-trivial deformation of $K$ that shows that this set is not rigid.

The boundary of $K$ in $\D$ consists of open  arcs of  circles
$C_i\sub \OC$, $i\in I$, that are orthogonal to $\partial \D$ and bound pairwise disjoint open disks 
$D_i\sub \OC$  in the complement of $K$. Here $I$ is some non-empty index set. Then $S= \OC\setminus\bigcup_{i\in I}D_i$ is a Schottky set in $\OC$
containing $K$. If we denote the reflection in the unit circle
$\partial \D$ by $R$, then $S=K\cup \geo K\cup R(K)$. For $i\in I$ denote by $R_i$ the reflection in $C_i$. Since $C_i$ is orthogonal to the unit circle, we have $R\circ R_i=R_i\circ R$. Let $\Gamma$ be the group generated by 
$R$ and $R_i$, $i\in I$. Then $\Gamma$ contains the Schottky group $\Gamma_S$ associated with $S$ as a subgroup of index $2$. Moreover, it follows that $$ S_\infty=\bigcup_{U\in \Gamma}U(Q)$$
is a measurable partition of $S_\infty$. Therefore, one can find a Beltrami coefficient $\mu$ on $\OC$ that is supported on $S_\infty$, that is invariant under $\Gamma$, and such that $\mu=\nu$ almost everywhere  on $\inte(Q)$ (cf.~the proof of Theorem~\ref{T:Rig2}). 
Let $F\:\OC\ra \OC$ be an orientation-preserving quasiconformal map with
$\mu_F=\mu$ almost everywhere .  By Lemma~\ref{L:gpinv} the map $F$ conjugates
$\Gamma$ to another group of M\"obius transformations. As in the proof of Lemma~\ref{L:SchtoSch} one sees that $R'=F\circ R\circ F^{-1}$ is a reflection in a circle. In particular, $F$ maps $\D$ to a disk. By post-composing the map $F$ by a M\"obius transformation if necessary (which does not change its Beltrami coefficient), we may assume that this disk is the unit disk. Then $F(\D)=\D$ and $F\circ R=R\circ F$. By 
Lemma~\ref{L:SchtoSch} the set $S'=F(S)$ is a Schottky set. Since $S$ is bounded by circles orthogonal to $\partial\D$ and $F$ commutes with $R$, the Schottky set $S'$ is also bounded by circles orthogonal to 
$\partial \D$. This implies that we can write $S'=K'\cup \geo K'\cup R(K')$,
where $K'=F(K)\in \Ca_2$.

As a quasiconformal map of $\D$ onto itself, the map $F$ is a quasi-isometry in the hyperbolic metric (this follows from standard distortion estimates for quasiconformal maps; see \cite[Ch.~9]{BHK} for more background).
 Since it maps  $K$ to another set in 
$\Ca_2$, it will follow that $K$ is not rigid, if we can show that there is no 
hyperbolic isometry on $\D$ that has finite distance to $F$ on $K$.
To see this we argue by contradiction, and suppose that there exists such an 
isometry. Then by Proposition~\ref{prop:qisom} there  exists a M\"obius transformation $\phi$ that leaves $\D$ invariant such that $F|\geo K=
\phi|\geo K$. Replacing $F$ by $\phi^{-1}\circ F$ if necessary, we may assume that $F$ is the identity on $\geo K$. Then $\geo K'=F(\geo K)=
\geo K$.  Since $K$ and $K'$ are the hyperbolic convex hulls of their boundaries at infinity, it follows that $K'=K$. So $F$ maps $K$ onto itself,
and is the identity on $\geo K$. But then $F$ is also a quasiconformal map of the quadrilateral $Q$ onto itself. Hence $\Mod(F(Q))=\Mod(Q)$.
On the other hand, $\mu_F=\nu$ almost everywhere on $\inte(Q)$, and so $\Mod(F(Q))\ne \Mod(Q)$ according to the choice of $\nu$. This contradiction shows that $K$ is not rigid.
\hfill $\Box$

\begin{theorem}\label{T:Rig5}
A set $K\in \Ca_3$ is rigid if and only if $\geo K$ has measure zero.
\end{theorem}

This statement corresponds to 
Theorem~\ref{T:Rig2}. 

\medskip
\no {\em Proof:} Let $K\in \Ca_3$, and $S=\geo K$. If $S$ has measure zero, 
Theorem~\ref{T:Rig3} implies that $K$ is rigid. 

Suppose $S$ has positive measure. In the proof of Theorem~\ref{T:Rig2} it was shown that there exists a quasiconformal map $f\: \Sph^2\ra \Sph^2$ such that 
$S'=f(S)$ is a Schottky set
 and such that $f|S\ne \ga|S$ for all M\"obius transformations $\ga$ on $\Sph^2$.  

    By the version of the Tukia-V\"ais\"al\"a theorem given in 
Proposition~\ref{prop:ballext}, there exists a quasisymmetric  map $F$ on the closed unit ball 
 extending $f$. The map $F$ is a quasi-isometry 
 on the open unit ball equipped with the hyperbolic metric. In this way, we obtain a quasi-isometry $F\: \H^3\ra \H^3$ with boundary map $\tilde F=f$.

 Let $K'\in \Ca_3$ be the hyperbolic convex hull of the Schottky set $S'$. 
 We claim that $F(K)$ and $K'$ have finite  Hausdorff distance. To see this, let $C$ and $C'$ 
 be the union of all geodesics with endpoints in 
 $S$ and $S'$, respectively. 
  Then $ \dist_H(K, C)<\infty$ and 
  $ \dist_H(K', C')<\infty$ (cf.~\cite[Proposition 10.1]{BS}). By geodesic stability of Gromov hyperbolic spaces, we also have $ \dist_H(F(C), C')<\infty$, and so $\dist_H(F(K), K')<\infty.$
 This implies that we can move each point in $F(K)$
   by a  bounded amount to a point in $K'$.
   In this way, we obtain a quasi-isometry $G\: K\ra K'$ with finite distance to  $F$. In particular, for the induced boundary map $\tilde G\: \geo K=S \ra \geo K'=S'$ we have $\tilde G=f|S$. So $G$ does not have finite  distance to the restriction of any isometry of $\H^3$ to $K$, because otherwise $\tilde G=f|S$ would agree with the restriction of a 
   M\"obius transformation to $S$.  
It follows that $K$ is not rigid. 
\hfill $\Box$ \medskip

\end{document}